\documentclass[%
 aip,
 amsmath,amssymb,
 reprint,%
]{revtex4-1}

\usepackage{graphicx}
\usepackage{dcolumn}
\usepackage{bm}
\usepackage{siunitx}
\usepackage{natbib}
\usepackage{comment}

\usepackage[utf8]{inputenc}
\usepackage[T1]{fontenc}
\usepackage{mathptmx}
\usepackage{etoolbox}
\usepackage{subcaption}
\usepackage{graphicx}

\usepackage{float} 
\usepackage{xcolor}

\definecolor{darkyellow}{RGB}{204, 204, 0} 

\usepackage{nomencl}
\usepackage{ifthen}
\renewcommand{\nomgroup}[1]{%
\ifthenelse{\equal{#1}{C}}{\item[\textbf{Constants}]}{%
\ifthenelse{\equal{#1}{G}}{\item[\textbf{Greek Letters}]}{%
\ifthenelse{\equal{#1}{S}}{\item[\textbf{Subscripts}]}{}}}
}

\makeatletter
\def\@email#1#2{%
 \endgroup
 \patchcmd{\titleblock@produce}
  {\frontmatter@RRAPformat}
  {\frontmatter@RRAPformat{\produce@RRAP{*#1\href{mailto:#2}{#2}}}\frontmatter@RRAPformat}
  {}{}
}%
\makeatother

\usepackage[switch]{lineno}

\begin{document}

\preprint{AIP/123-QED}


\title[Modal Analysis of the Wake Shed Behind a Horizontal Axis Wind Turbine with Flexible Blades]{Modal Analysis of the Wake Shed Behind a Horizontal Axis Wind Turbine with Flexible Blades}

\author{S. Salavatidezfouli}
\altaffiliation[ ]{These authors contributed equally to this work.}
\affiliation{
Mathematics Area, MathLab, International School for Advanced Studies (SISSA), Trieste, Italy.}

\author{A. Sheidani}%
\altaffiliation[ ]{These authors contributed equally to this work.}
\affiliation{ 
Mathematics Area, MathLab, International School for Advanced Studies (SISSA), Trieste, Italy.
}%

\affiliation{ 
Department of Applied Physics, Eindhoven University of Technology, Netherlands.}

\author{K Bakhshaei}
\altaffiliation[ ]{These authors contributed equally to this work.}
\affiliation{ 
Mathematics Area, MathLab, International School for Advanced Studies (SISSA), Trieste, Italy.}

\author{A Safari}
\affiliation{%
Department of Mechanical Engineering, The Hong Kong University of Science and Technology, Hong Kong, China
}%

\author{A Hajisharifi}
\affiliation{ 
Mathematics Area, MathLab, International School for Advanced Studies (SISSA), Trieste, Italy.}

\author{G Stabile}
\email{giovanni.stabile@uniurb.it}
\affiliation{ 
Department of Pure and Applied Sciences, Informatics and Mathematics Section, University of Urbino Carlo Bo, Urbino, Italy}
\affiliation{ 
The Biorobotics Institute, Sant'Anna School of Advacended Studies,
Pisa, Italy}

\author{G Rozza}
\affiliation{ 
Mathematics Area, MathLab, International School for Advanced Studies (SISSA), Trieste, Italy.}

\date{\today}

\begin{abstract}
The proper orthogonal decomposition has been applied on a full-scale horizontal-axis wind turbine to shed light on the wake characteristics behind the wind turbine. In reality, the blade tip experiences high deflections even at the rated conditions which definitely alter the wake flow field, and in the case of a wind farm, may complicate the inlet conditions of the downstream wind turbine. The turbine under consideration is the full-scale model of the National Renewable Energy Laboratory 5MW onshore wind turbine which is accompanied by several simulation complexities including turbulence, mesh motion and fluid-structure interaction. Results indicated an almost similar modal behaviour for the rigid and flexible turbines at the wake region. In addition, more flow structures in terms of local vortices and fluctuating velocity fields take place at the far wake region. The flow structures due to the wake shed from the tower tend to move towards the center and merge with that of the nacelle leading to an integral vortical structure 2.5 diameter away from the rotor. Also, it is concluded that the exclusion of the tower leads to missing a major part of the wake structures, especially at far-wake positions. 
\end{abstract}

\maketitle

\section{Introduction}


Clean energies are considered to be a suitable substitute for conventional carbon-intensive energy resources to decrease environmental disasters such as climate change, global warming and acid rain. In this matter, a great portion of the global energy demand due to the world population surge in the last decades can be addressed with renewable energy resources \citep{alkhabbaz2022impact, Salavatidezfouli2024}. Therefore, more investments have been drawn to the development and optimization of utilization of these energy sources \citep{bajuri2022computational}. 

Wind energy is characterized as a feasible substitute owing to its availability, reliability and cleanness. Wind Turbines are one of the major sources of energy production burdening low impact on the environment and are associated with high efficiency and performance\citep{MAJIDINEZHAD2022112791}. This led to a great rise in the installation of wind turbines around the world and it is predicted that the share of wind turbines in electricity production will rise to 5\% by 2030 \citep{moghaddam2011wind}. To meet the increasing demand for energy production, large-scale wind turbines are going to be clustered in more extensive areas with condensed configurations to ensure enough energy production. In fact, two challenges are always highlighted among the vital design factors of the large-scale offshore wind farm \citep{uchida2021doppler}: (i) the wake shed by the upstream turbines which can decrease the amount of power generated by the downstream turbines and (ii) the probability of failure of the wind turbine structures including nacelle and blades. An efficient and reliable design of wind turbine heavily depends on the identification of the dynamics of wind turbine wakes and their effects associated with neighbouring wind turbines for optimization and development of next-gen wind turbines and wind farms \citep{uchida2022effects}. Moreover, a full analysis of the near and far wake regions of the wind turbine is necessary to accomplish the performance assessment of the blades and further acoustic analysis. The wake region behind a wind turbine can exhibit flow structures ranging from millimetres to hundred meters which correspond to the size of wind turbine components \citep{purohit2021accuracy}. 

Although large-scale offshore wind turbines undergo severe aerodynamic and weight forces, they have attained a lot of attention during recent years due to their high ability in power production, less visual and environmental impacts and abundance of wind in offshore areas \citep{liu2019aeroelastic}. To overcome this complication, manufacturers aim to produce lighter blades to mitigate the problem of large loads applied to the tower and mechanical connections \citep{liu2022finite}. As a result, the larger and lighter blades lead to more flexible structures which alter the performance of the wind turbine and render the assumption of rigid blades invalid. \citep{hand2021structural}.

\textcolor{darkyellow}{Due to spatial limitations, the future trajectory of wind farms involves deploying wind turbines in close proximity within a tandem configuration. Therefore, achieving precise modelling of wake dynamics becomes imperative to optimize layout, operation, and control. The intricate interplay between wake dynamics and rotor aeroelasticity significantly influences each other. Conventional aerodynamic theories fall short in adequately capturing the interaction between the fluctuating inflow and the dynamic aerodynamics of the flexible blades \citep{ma2019analysis}. Blades traversing through the air shed vorticity, typically discharged downstream and away from a rigid structure. However, when the blades deform out-of-plane, the rotor engages with its vorticity, casting doubt on the accuracy of design assumptions. Moreover, the structural dynamics of blades, inherent curvature and sweep, and significant nonlinear deflection (encompassing torsion and bend-twist coupling) further complicate the physics and the evaluation of stability \citep{della2022two}.
Consequently, this study aims to investigate the wake shed phenomenon behind a flexible horizontal-axis wind turbine, with Proper Orthogonal Decomposition (POD) serving to extract hidden fluid flow structures.}

The employment of POD to extract different modes of flow dynamics was first proposed by Lumely et al. \citep{lumley1967structure}. This idea was then nurtured and developed into more robust methods in different studies namely, Sirovich \citep{Sirovich_1987} who put forward the concept of snapshot POD. Consequently, in recent years, the advances in Reduced Order Modelling (ROM) techniques have led to plenty of studies intended to enlighten the hidden flow characteristics by employing methods like POD \citep{benner2021system, hajisharifi2023non, hajisharifi2024comparison, hajisharifi2024lstm, salavatidezfouli2024predictive, StarStabileRozzaDegroote2020}. For instance, Premaratne et al. \citep{premaratne2022proper} investigated coherent vortices embedded in the turbulent wake flow behind a prototype Horizontal Axis Wind Turbine (HAWT) experimentally. A high-resolution Particle Imaging Velocimetry (PIV) system was used to quantify the evolution of the unsteady wake vortices. The POD analysis was followingly conducted to identify energetic vortex patterns in the turbulent near wake. They concluded that the wake region adjacent to the turbine accounts for a high percentage (62\%) of the modal energy. While highly stretched fluctuations, especially in the vertical direction, were observed at the end of the wake region. 

Due to the large dimensions of such turbines, the experimental study regarding the blades' flexibility effects is impractical. In addition, it is hard to establish a complete similarity between the actual model and the prototype. For instance, in the existence of sea waves, an equal value of Froude number should be maintained to ensure hydrodynamic similarity. However, maintaining similar aerodynamic dimensionless numbers such as the Tip Speed Ratio (TSR) and Reynolds Number at the same time is almost impossible \citep{Fontanella_2021, Madsen_2020}. As a result, the numerical methods may yield more reliable outcomes for offshore wind turbines owing to the fact that all the details are taken into consideration. 

\textcolor{darkyellow}{In recent years there has been an increasing interest in studies aiming at understanding the advanced working conditions of the wake of the flexible wind turbine. To illustrate, Zhang et al.\cite{ZHANG2024105625} proposed a novel two-way fluid-structure interaction (FSI) model for the simulation of a flexible HAWT taking the effect of inflow condition into consideration as well. As one of the main findings, it was demonstrated that the tower has a profound effect on the wake of the wind turbine causing the vortical behaviour of the flexible blades to undergo considerable alterations. Farrell et al. \citep{farrell2024analyzing} performed a numerical simulation to study the wake of a wind farm of flexible HAWTs under atmospheric turbulent oscillations conditions. The results proved that the wake affects the vibration dynamics of the wind turbines located downstream significantly. Trigaux et al. \cite{trigaux2024investigation} employed the flexible actuator line method implementing the Large Eddy Simulation (LES) technique to investigate the aeroelasticity of a 15MW offshore HAWT. It was reported that in general, an assumption of a rigid blade is sufficient to capture the wake dynamics. Yu et al.\cite{yu2020aeroelastic} utilized the elastic actuator line method which regards the blade as an elastic beam for the FSI application to study the performance and wake of two wind turbines in an aligned configuration. It was shown that the flexibility of the blades along with the wake of the tower can lead to instability occurrence on the downstream turbine. Liu et al.\cite{liu2019aeroelastic} studied a flexible offshore wind turbine with platform surge-induced motion. The results proved that the blade flexibility without considering the effect of platform motion does not change the performance of the wind turbine more than $5\%$. Leng et al.\cite{leng2023fluid} employed the actuator line method to perform a two-way FSI study of large HAWTs. It was reported that blade flexibility can help the wake recovery process compared to rigid cases.}

From the numerical perspective, to take the effect of blade flexibility into consideration, the structural and fluid dynamics phenomena need to be modelled by independent solvers followed by an interface to transmit the results from one to the other \citep{deng2020investigating}. Obviously, the accuracy and computational cost are the two main factors which govern the FSI simulations, and one needs to strike a balance to conduct a cost-efficient study. For instance, Blade Element Momentum (BEM) oriented techniques are mainly employed to reduce the cost of the computational cost associated with the flow simulation \citep{Rafiee_2016}. The FAST code was first developed by the National Renewable Energy Laboratory (NREL) to extend the aerodynamic calculation of horizontal off-shore wind \citep{prowell2010fast}. The code utilized a variety of packages to address fluid and structure models, and more inclusively BEM method for simulating the fluid flow. Several research studies were reported on the utilization of the FAST package to model FSI in wind turbines \citep{borouji2019fluid}. 

\textcolor{darkyellow}{Due to the simplifications related to the BEM theory, many attempts have been made to modify and improve the accuracy of such simulations. To illustrate, Mo et al. \citep{Mo_2015}, and Rajan et al.\citep{Rajan_2018} modified the BEM theory by considering the fluxional fluctuations and torsional deformations for each element. Yu and Kwon employed a non-linear Euler–Bernoulli beam model coupled with an in-house incompressible flow solver in their analysis \citep{yu2014time}. They demonstrated that the aeroelastic deformation of NREL 5MW wind turbine, boasting a rotor diameter (D) of $126m$, decreases aerodynamic loads through torsional deformation. Rajan et al. \citep{Rajan_2018} developed the Dynamic Rotor Dynamics (DRD) – BEM models to capture the three-dimensional effects of blade deformation more accurately by constantly altering the point of element equilibrium. Sayed et al. \citep{Sayed_2019} conducted a numerical study on a DTU 10MW HAWT by combining Unsteady Reynolds Average Navier-Stokes (URANS) with a structural solver called Carat++ which calculates the deformation of the blade a non-linear-Timoshenko-beam assumption. In this work, an increase in the wind turbine power due to the blade flexibility by 1\% was reported. On the other hand, the discrepancy between the rigid and flexible simulations has been reported to be more significant for other HAWT types. For instance, Badshash et al. \citep{Badshah_2018} showed almost a 5\% difference in the power coefficient of a two-bladed HAWT by considering the effect of blade flexibility. Generally, the effect of blades' flexibility plays an important role in larger turbines.}

The precise perception of flow characteristics in the wake region of the wind turbine requires the employment of accurate numerical methods accompanied by powerful analysis tools. \textcolor{darkyellow}{Even though cost-efficient methods like BEM or simplified structure concluded high accuracy in terms of output power, they are not able to accurately predict the wake region behind the wind turbine. Hence, as the main objective of this paper, we aim to implement a fully coupled high-fidelity FSI simulation to account for the wake characteristics behind a large-scale HAWT. On the other hand, POD analysis was implemented to reveal flow structures in the wake region as it is claimed to be a great tool in this respect  \citep{sheidani2023assessment}. Additionally, as another objective, to understand to what extent the blade's flexibility alternates wake properties, a comparison between the rigid and flexible blades was carried out. NREL 5MW horizontal wind turbine was selected for these goals based on the work of Dose et al. \citep{dose2020fluid}. At first, a validation study for the FSI solver with respect to the numerical results of Dose et al.\citep{dose2018fluid} is performed. Afterwards, the POD analysis of the wake region based on the high-fidelity simulation of rigid and flexible turbines is presented to unfold the physics of the wake followed by a conclusion on the results.}

\section{\label{sec:level2}Numerical Method}

The rigid and flexible simulations of the wind turbine have been carried out with STAR CCM+ solver. It incorporates \textit{Finite Volume Method (FVM)} and \textit{Finite Element Method (FEM)} solvers, respectively, for the flow and structural simulations. In addition, mesh motion, overset and mesh morphing techniques were employed simultaneously to address the nodes' movement through the simulation. In the following section, the details of the utilized solvers will be discussed.

\subsection{Solver Description}
\paragraph{Fluid Flow.}%
The incompressible fluid flow assumption is often utilized for wind turbines. Especially, in the case of large-scale horizontal wind turbines, the Mach number is typically low which assures incompressibility conditions. In order to perform the numerical simulation of the problem, the URANS was numerically solved.

The URANS formulation adopted for the continuity and momentum equations of the incompressible fluid is as follows  \citep{SALAVATIDEZFOULI2024Effect}:
\begin{equation}
\label{eq1}
\begin{aligned}
&\frac{\partial \overline{U_{i}}}{\partial x_i}=0,\\
\end{aligned}
\end{equation}
\begin{equation}
\label{eq2}
\begin{aligned}
&\rho \frac{\partial \overline{U_i}}{\partial t}+\rho \frac{\partial \overline{U_i} \; \overline{U_j}}{\partial x_j}=\\
&\;\;\;\;\;\;\;\;-\frac{\partial P}{\partial x_i}+\frac{\partial}{\partial x_j}\left[\mu\left(\frac{\partial \overline{U_i}}{\partial x_j}+\frac{\partial \overline{U_j}}{\partial x_i}\right)-\rho \overline{u_{i}^{\prime} u_{j}^{\prime}}\right],\\
\end{aligned}
\end{equation}
where $\overline{U_{i}}$ and $\rho \overline{u_{i}^{\prime} u_{j}^{\prime}}$ are the averaged velocity and Reynolds stress tensor, respectively. The latter can be written as:
\begin{equation}
\label{eq3}
\begin{aligned}
&-\rho \overline{\boldsymbol{u}_{i}^{\prime} \boldsymbol{u}_{\boldsymbol{j}}^{\prime}}=2 \mu_t \overline{S_{ij}}-\frac{2}{3} k \rho \overline{\delta_{ij}},\\
\end{aligned}
\end{equation}
where the mean strain rate tensor is:
\begin{equation}
\label{eq4}
\begin{aligned}
&\overline{S_{ij}}=\frac{1}{2}\left(\frac{\partial \bar{U}_i}{\partial x_j}+\frac{\partial \bar{U}_j}{\partial x_i}\right).\\
\end{aligned}
\end{equation}
The turbulent viscosity term, $\mu_t$, can be addressed with common one- or two- equation eddy viscosity models to provide turbulence closure \citep{safari2023numerical, sheidani2022study}.  

The selection of the suitable Reynolds Average Navier-Stokes (RANS) model is vital for HAWT as separation occurs during the blades' rotation. Many reports have agreed on both the accuracy and simplicity of $k-\omega$ Shear Stress Transport (SST) model for capturing turbulent features of the flow around the wind turbines \citep{cai2022effects}, and hence, was utilized for this research. All the spatial and temporal terms in Equation \ref{eq2} were discretized by the second-order method and the Semi-Implicit Method for Pressure-Linked Equations (SIMPLE) algorithm was employed to solve the set of equations, as suggested by literature \citep{salavatidezfouli2023deep, bakhshaei2021multi}. 

\paragraph{Structural.}%
The deformation of the flexible blades is described using a Lagrangian formulation and is governed by the structural momentum conservation equation, which can be expressed as follows:
\begin{equation}
\label{eq_solid1}
M_s \ddot{d}+C_s \dot{d}+K_s d=F_{s}, \\
\end{equation}
where $M_s$, $C_s$ and $K_s$  represent the mass, damping coefficient and stiffness of the turbine blade, respectively. $\ddot{d}$, $\dot{d}$ and $d$ are acceleration, velocity, and displacement of the structure element, respectively.

$F_s$ represents the force vector acting on the solid including pressure, shear, Coriolis, centrifugal and mass forces:
\begin{equation}
\label{eq_solid2}
F_{s}=F_{P} + F_{\tau} + F_{Cori.} + F_{cent} + F_{w}. \\
\end{equation}
The solution was derived through time integration by solving the motion equation in the Newmark scheme and the solution to the solid displacement was formulated in terms of the time history. More details on the implementation of fluid-structure interaction can be found in Appendix \ref{two-way coupling}.

\subsection{Model Description}
\begin{figure}
\nolinenumbers
	\centering
	\includegraphics[scale=0.09]{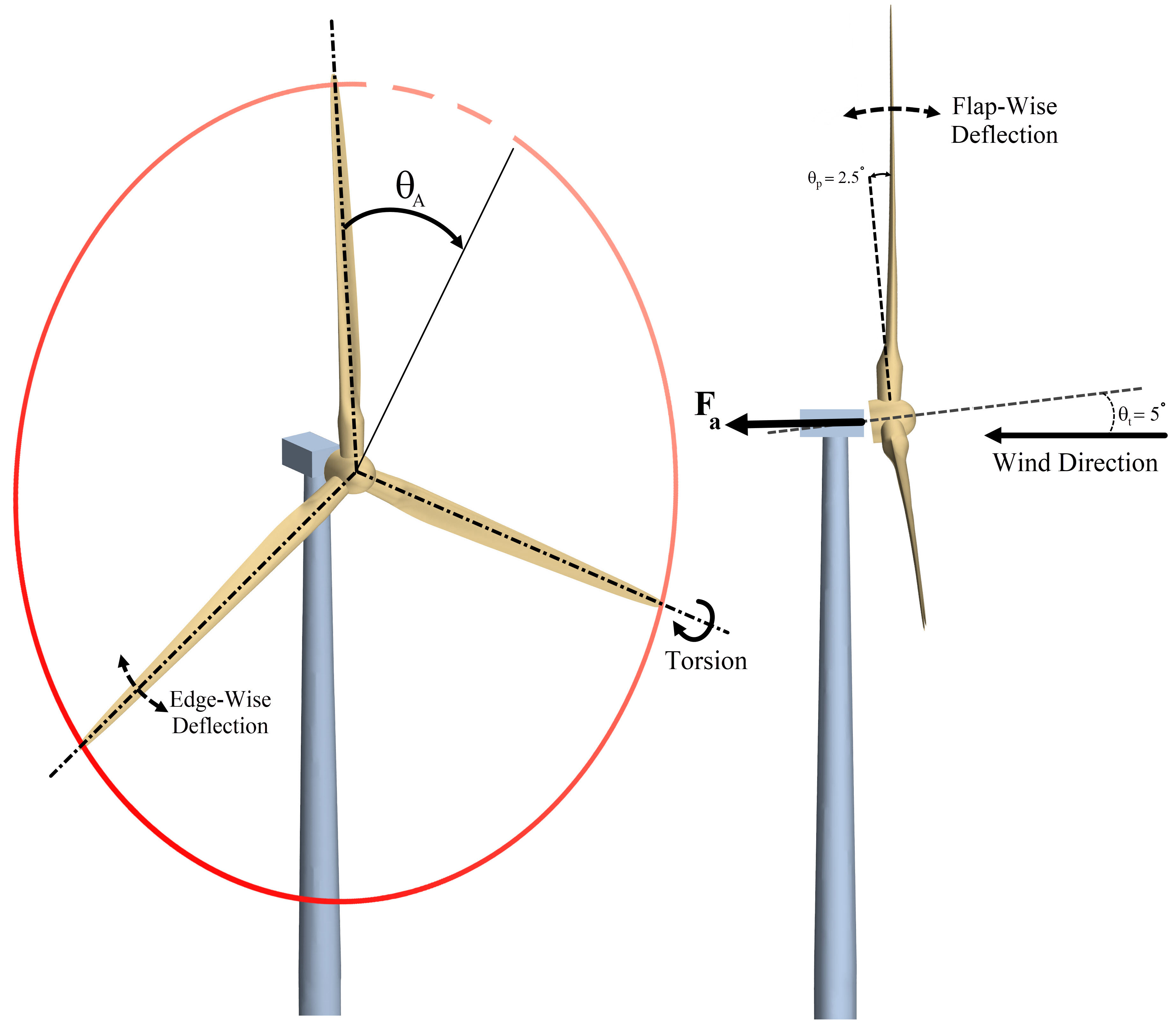}
	\caption{Geometrical properties of the NREL 5MW}
	\label{Fig_Blade}
\end{figure}

The wind turbine model used in the present work is the onshore NREL 5 MW based on the work of Dose et al. \citep{dose2020fluid}, as shown in Fig. \ref{Fig_Blade}. The detailed specifications of the turbine model are presented in Table \ref{tab:Specifications}. Other geometrical details of the turbine can also be found in the works of Coulling et al. \citep{coulling2013validation} and Zhang et al. \citep{zhang2018fully}.
\begin{table}[]
\nolinenumbers
\centering
\renewcommand{\arraystretch}{1.5}
\caption{Specifications of the NREL 5MW onshore wind turbine}
\label{tab:Specifications}
\begin{tabular}{|c|c|}
\hline
Rotor Orientation,   configuration & Upwind, 3-Blades     \\ \hline
Blade Sections                     & DU and NACA airfoils \\ \hline
Blade Length                       & 61.5m               \\ \hline
Rotor Diameter $(D)$                   & 126m                \\ \hline
Tilt Angle ($\theta_t$)                         & 5$^{\circ}$                    \\ \hline
Precone Angle ($\theta_P$)                      & 2.5$^{\circ}$                  \\ \hline
Hub Diameter                       & 3m                  \\ \hline
Tower Height                       & 90m                 \\ \hline
\end{tabular}
\end{table}

\begin{figure*}
\nolinenumbers
	\centering
	\includegraphics[scale=0.08]{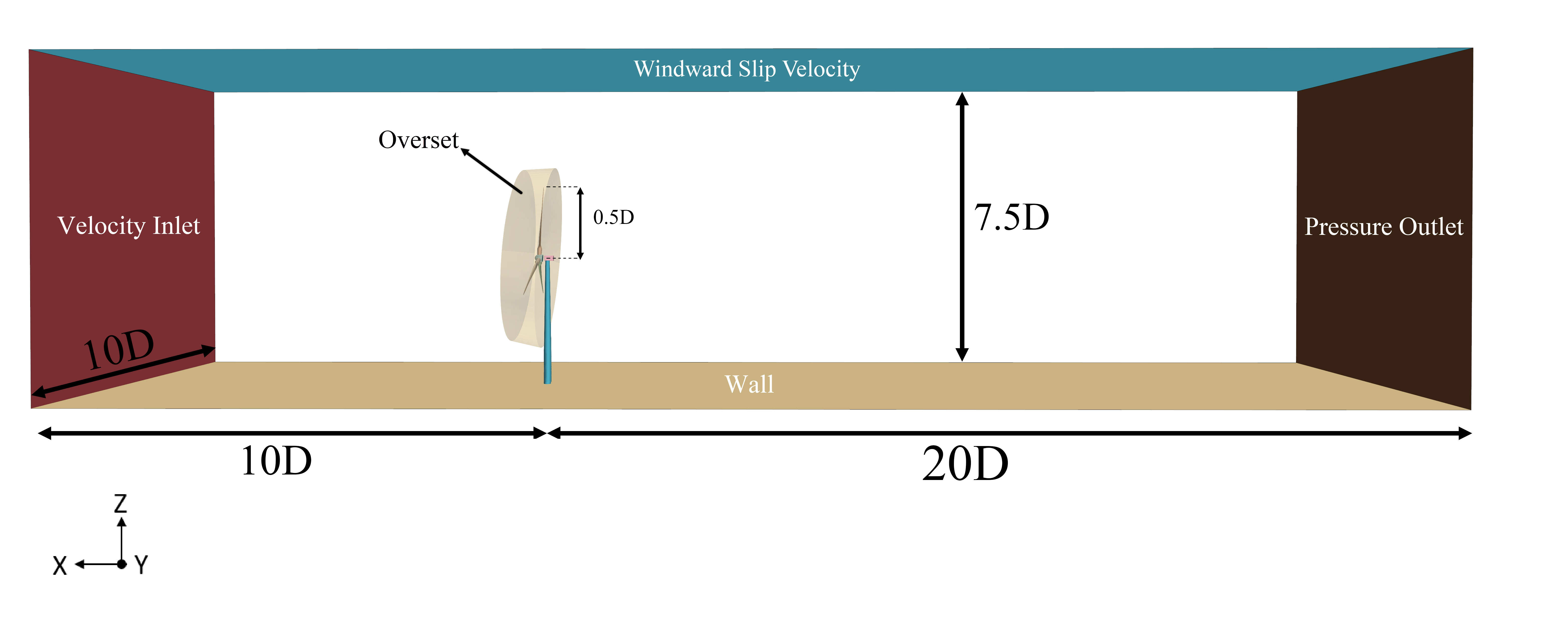}
	\caption{\textcolor{darkyellow}{Fluid domain and boundary conditions considered for the onshore NREL 5MW}}
	\label{FIG_Domain1}
\end{figure*}

Fig. \ref{FIG_Domain1} illustrates the computational domain of the fluid. The upstream and sides boundaries are defined as the velocity inlet boundary condition and the pressure outlet is considered for the downstream boundary. The no-slip wall condition is imposed on the turbine walls and the bottom boundary. All the boundaries were put far enough from the turbine to ensure flow uniformity. Thereby, a hexahedral domain with the size of $3000\si{m} \times 1000\si{m} \times  750\si{m}$ (length, width, height) was generated. 

The rated conditions presented by Jonkman et al. \citep{jonkman2009definition} were considered for the wind turbine. In this regard, a uniformly distributed wind velocity and rotor angular velocity of $11.4\si{m/s}$ and 12.1$\si{rpm}$ were employed, respectively. A turbulence intensity of 0.14 was considered at the inlet boundary. The air density and viscosity were set to $1.2 \, \si{kg/m^3}$ and $1.8\times {10^{-5}} \si{Pa.s}$, respectively.
%
%
%

It is worth mentioning that the current research study only focuses on the effect of blade deformation on the wake properties behind the wind turbine. Thus, there would be no need to model the deformation of the tower and nacelle. Structural properties were considered similar to that of Asareh et al. \citep{asareh2013evaluation}.

\subsection{Mesh Generation}
\begin{figure*}
\nolinenumbers
	\centering
	\includegraphics[scale=0.08]{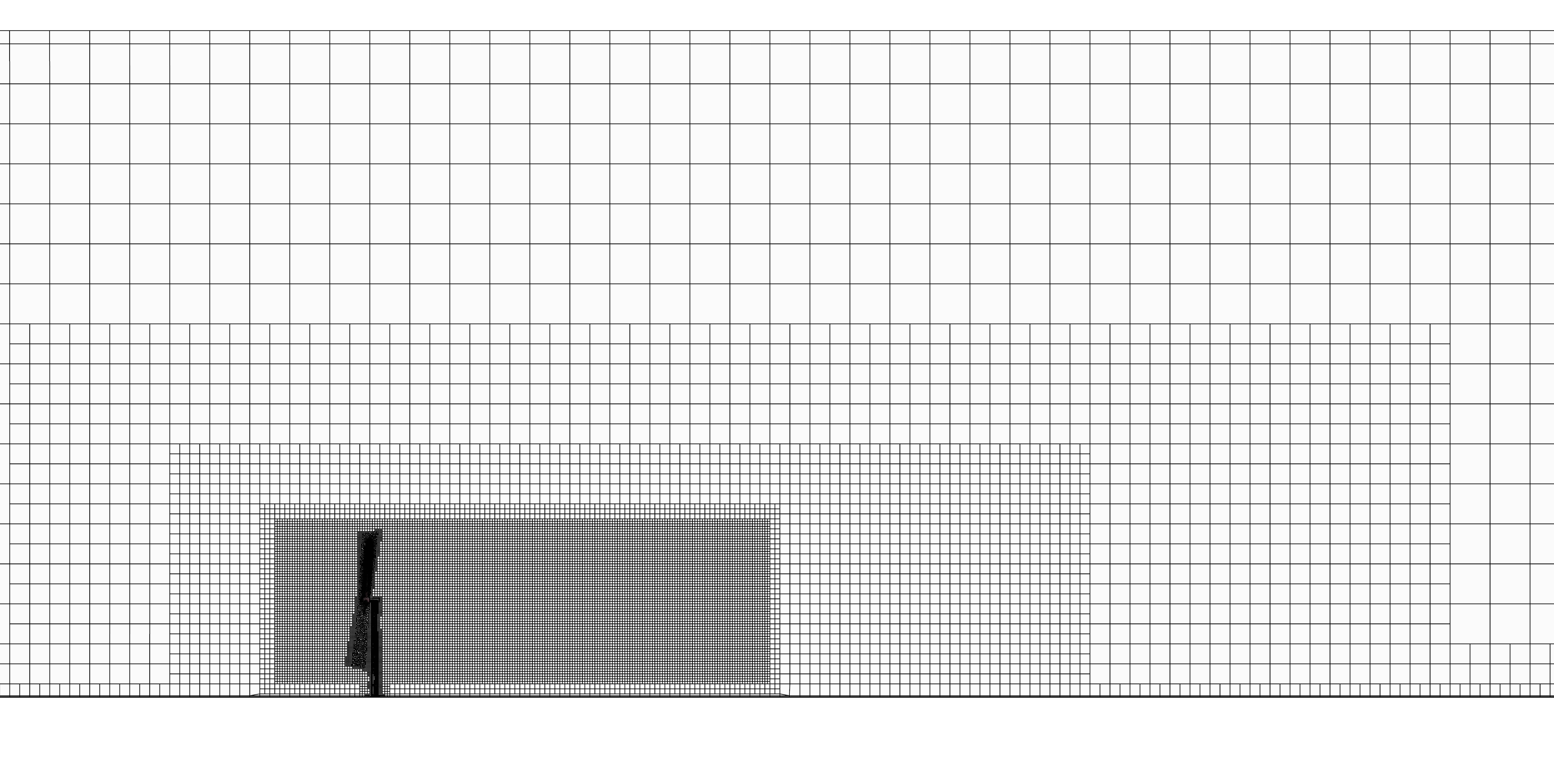}
	\includegraphics[scale=0.115]{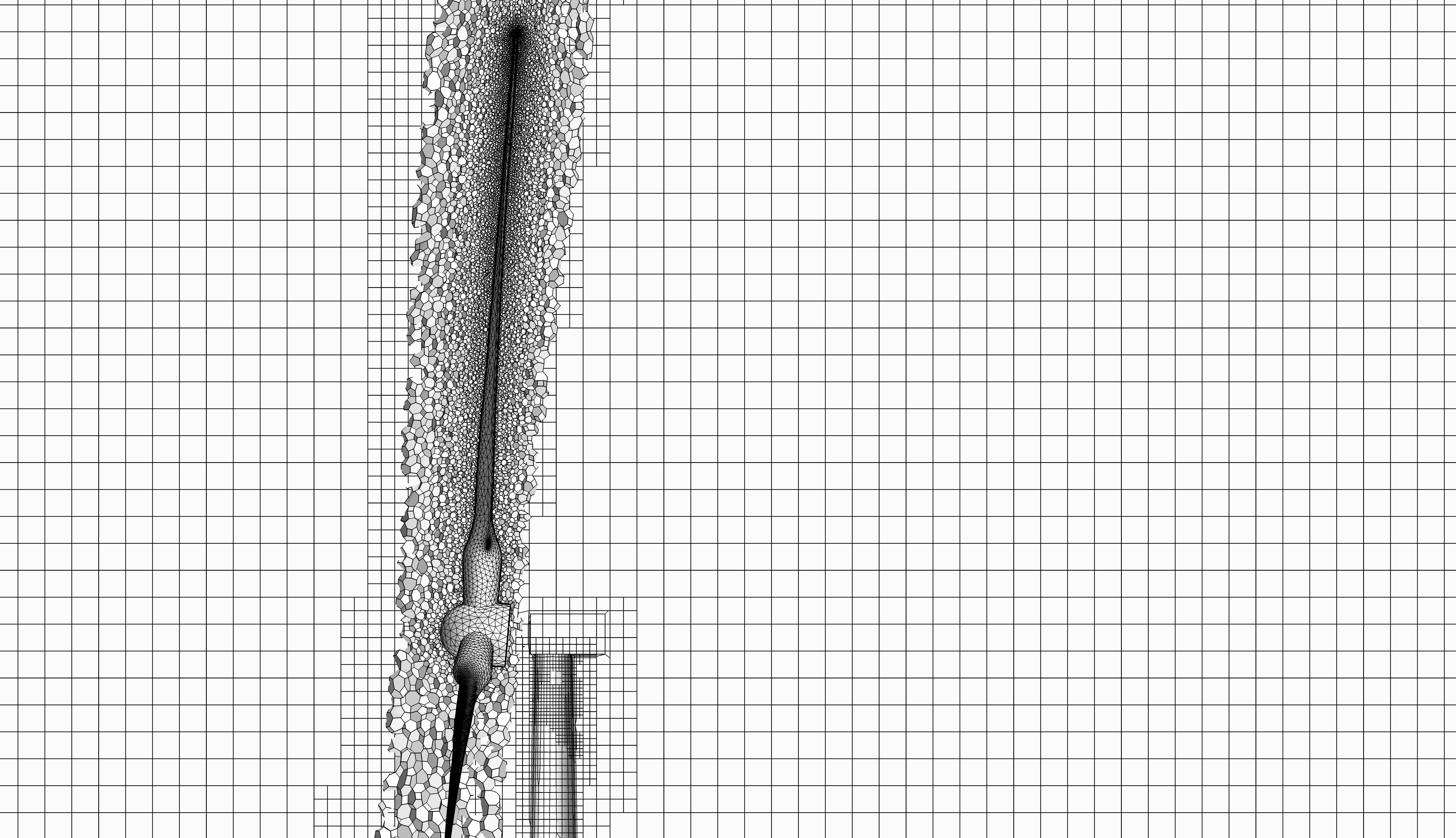}
	\caption{Mesh of the background and overset regions }
	\label{FIG_MeshDomain_Overset}
\end{figure*}
\begin{figure*}
\nolinenumbers
	\centering
	\includegraphics[scale=0.0725]{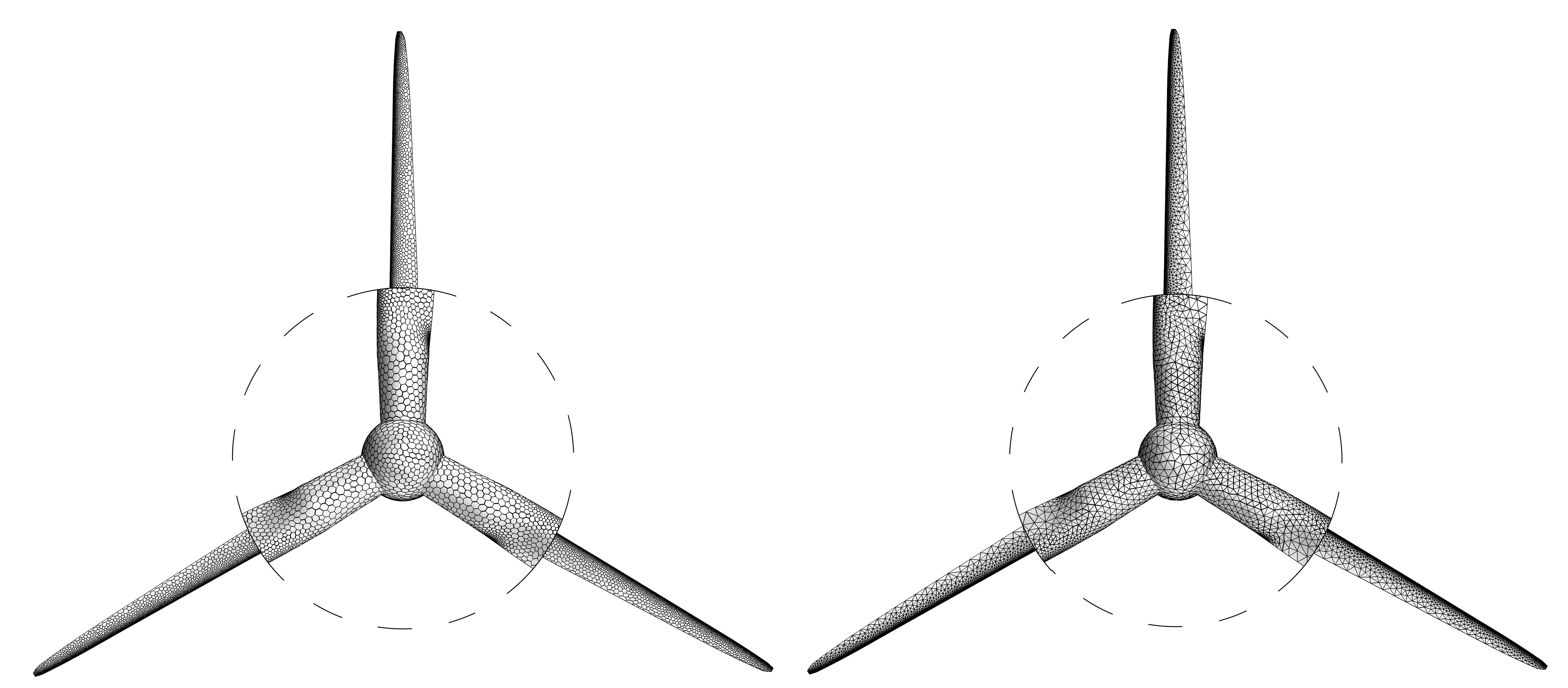}
	\caption{Polyhedral and tetrahedral meshes for the blade in the fluid and solid regions, respectively}
	\label{Fig_BladeMeshes}
\end{figure*}
Precise considerations need to be taken for the mesh generation of the large-scale blade-resolved simulation of wind turbines. While it is required to implement either the sliding mesh or overset technique for rigid wind turbines, flexible FSI simulations impose an additional complexity caused by the blade deformation. To address this difficulty, a combination of the overset technique with mesh morphing was implemented.

As for Computational Fluid Dynamics (CFD), hexahedral elements are suited for simplified geometries which can guarantee a better orthogonal quality, and therefore, result in better accuracy, whereas the polyhedral elements are suited for models with geometrical complexities such as the blade tip and leading edge \citep{wang2021comparison}. For implementing the overset solver, two fluid regions namely background and overset, as shown in Fig. \ref{FIG_MeshDomain_Overset}, were created. Hexahedral elements were utilized for the rectangular background and for the overset region, which encompasses the blades, polyhedral elements with a refined boundary layer mesh were considered. Fig. \ref{FIG_MeshDomain_Overset} shows the corresponding meshes for the fluid regions.

The solver supports hexahedral and tetrahedral elements. The latter was utilized for the structural model of the wind turbine, due to the mentioned reasons. A comparison between the fluid and structural meshes of the blade is depicted in Fig. \ref{Fig_BladeMeshes}.

\subsection{Verification and Validation}
The mesh sensitivity analysis needs to be carried out for the output parameters to asses meshing-related numerical error. To do so, coarse, medium and fine meshes with the total number of $2.4 \times 10^{6}$, $4.6 \times 10^{6}$ and $5.9 \times 10^{6}$ elements were created and examined. Table \ref{tab:meshInfo} presents a detailed breakdown of the number of elements for each individual region. 

As for the dimensionless wall distance index, $y^{+}$, the average value was kept below 40 for all cases, as instructed by Siddiqui et al. \citep{siddiqui2017quasi}, with at least 10 prism layers \citep{lillahulhaq2019numerical} to capture the near-wall phenomena and optimize the element number, simultaneously. 

\begin{figure}
\nolinenumbers
\centering
\includegraphics[scale=0.35]{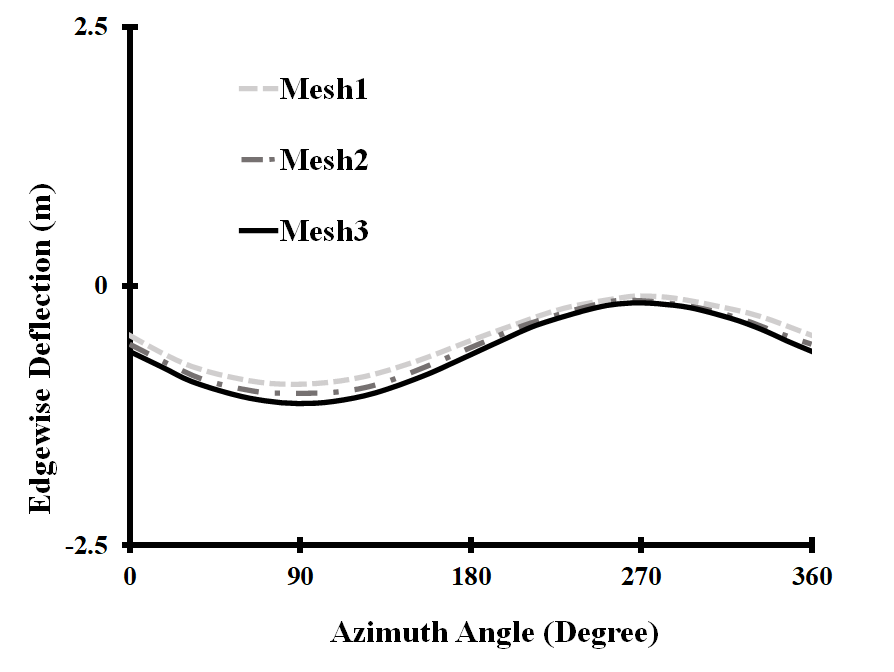}
\caption{Comparison of the edge-wise deflection in one cycle for three meshes}
\label{Fig_EdgeDef}
\end{figure}

A comparison between edgewise deflection for three meshes is shown in Fig. \ref{Fig_EdgeDef}. There exists a negligible difference between the results of the second and third mesh, and hence, the rest of the post-processing including POD extraction was carried out for the second mesh.
\begin{table}[]
\nolinenumbers
\centering
\renewcommand{\arraystretch}{1.5}
\caption{\textcolor{darkyellow}{Number of elements for different regions of the model}}
\label{tab:meshInfo}
\begin{tabular}{|c|c|c|c|}
\hline
\textbf{}        & \textbf{Mesh \#1} & \textbf{Mesh \#2} & \textbf{Mesh \#3}  \\ \hline
\textbf{Structure}            & 305915            & 593197            & 794005            \\ \hline
\textbf{Fluid Overset}    & 788520            & 1436501           & 2391163           \\ \hline
\textbf{Fluid Background} & 1362236           & 2561289           & 2752772           \\ \hline
\textbf{Total} & 2456671           & 4590987           & 5937940           \\ \hline
\end{tabular}
\end{table}

\begin{figure}[h]
\nolinenumbers
	\centering
	\includegraphics[scale=0.3]{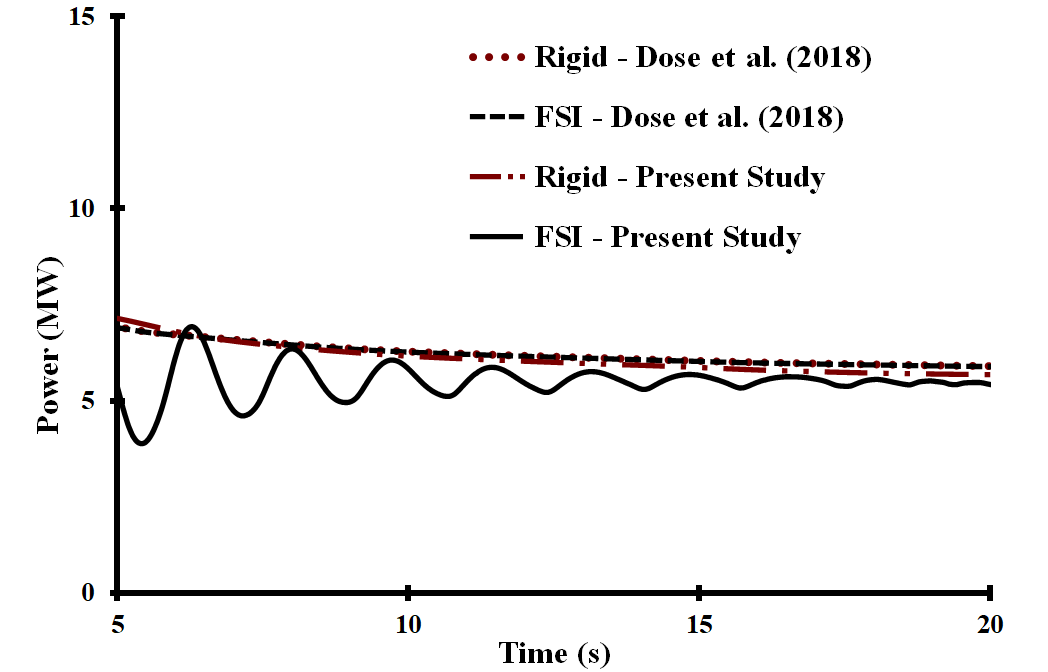}
	\caption{Evolution of the turbine power over the simulation time for rigid and FSI simulations}
	\label{Fig_Power}
\end{figure}

To validate the flow solver the generated power of the wind turbine for the rigid and FSI cases was compared with the results of Dose et al.  \citep{dose2018fluid}, which is shown in Fig. \ref{Fig_Power}. As mentioned earlier, the oscillatory nature of FSI results corresponds to the change of elastic properties of the blade at the initial times (see section \ref{two-way coupling}. After approximately three revolutions (i.e. $T = 15\si{sec}$), the numerical results almost reached that of Dose et al. \citep{dose2018fluid}.
\begin{figure}[h]
\nolinenumbers
	\centering
	\includegraphics[scale=0.25]{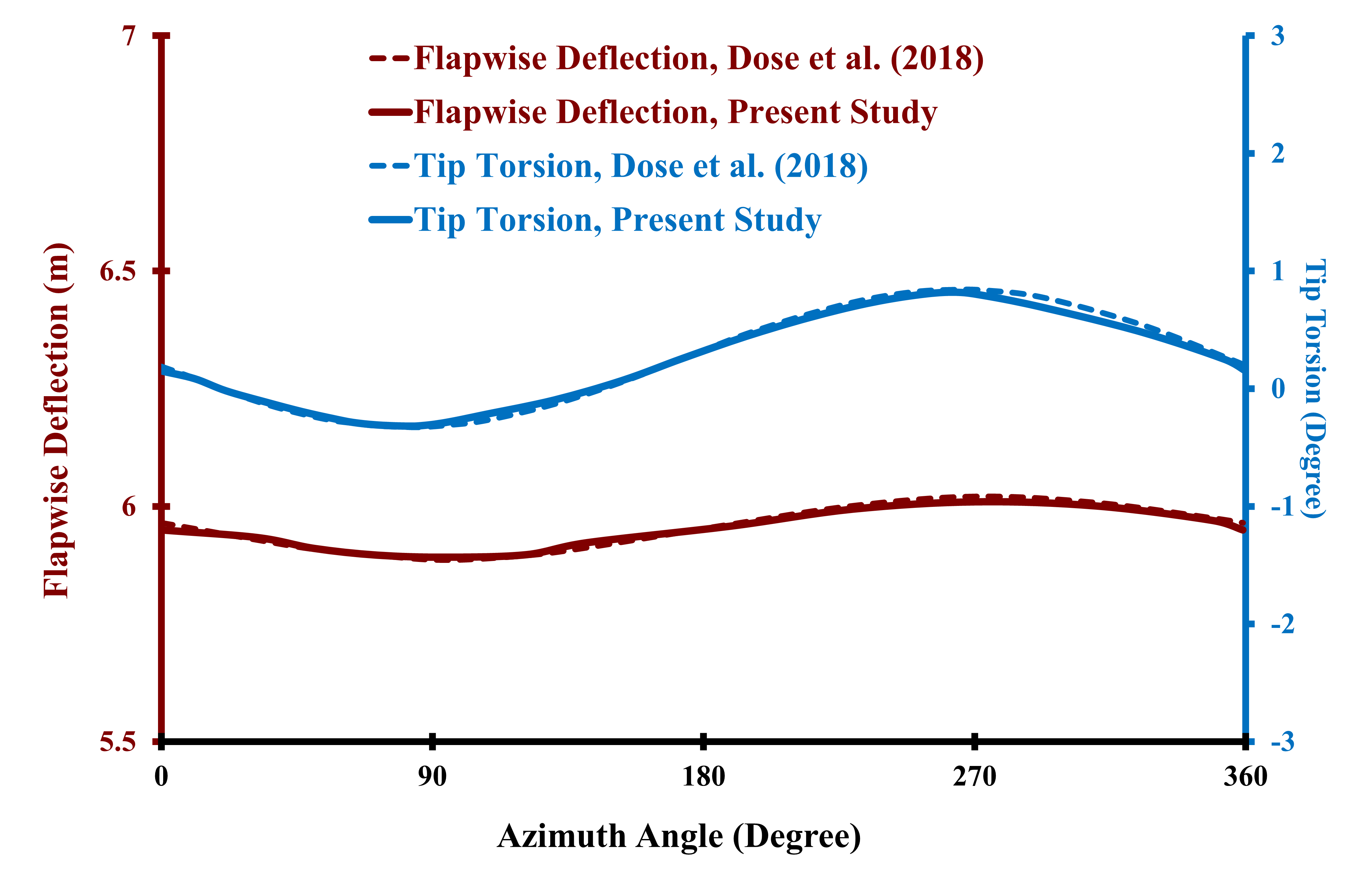}
	\caption{\textcolor{darkyellow}{The flap-wise deflection and tip torsion of a blade through one rotation cycle}}
	\label{Fig_FlapDefTipTors}
\end{figure}

The flap-wise deflection and tip torsion were selected to assess the FEM solver. The tip deflection and tip torsion through one revolution of the turbine are shown in Fig. \ref{Fig_FlapDefTipTors} and compared with those of Dose et al. \citep{dose2018fluid}. These conclude that the FSI solver could accurately capture the blade deformation through time. 

\begin{table*}[]
\nolinenumbers
\centering
\renewcommand{\arraystretch}{1.5}
\caption{\textcolor{darkyellow}{Comparison of the results under the rated working condition of NREL 5MW wind turbine}}
\label{tab:validation}
\nolinenumbers
\begin{tabular}{|c|c|c|c|c|}
\hline
                            & \textbf{Edge-Wise} & \textbf{Rotor} & \textbf{Thrust} \\ 
                         & \textbf{ Deflection (m)} &\textbf{Power (MW)} &\textbf{Force (kN)}\\\hline
\textbf{Sabale and Gopal \citep{sabale2019nonlinear}}                          & -0.57                             & 4.97                      & 490                        \\ \hline

\textbf{Dose et al. \citep{dose2018fluid}}                         & -0.62                             & 5.49                      & 771                        \\ \hline
\textbf{Jonkman et al. \citep{jonkman2009definition}}                         & -0.61                             & 5.28                      & 814                        \\ \hline
\textbf{Yu and Kwon \citep{yu2014predicting}}                           & -0.63                             & 5.22                      & 656                        \\ \hline
\textbf{Ponta et al. \citep{ponta2016effects}}                           & -0.56                             & 5.19                      & 660                        \\ \hline
\textbf{Li et al \citep{li2020aerodynamic}}                            & -0.57                             & 5.3                       & 678                        \\ \hline
\textbf{Present Study}  & -0.58                             & 5.01                      & 636                        \\ \hline
\end{tabular}
\end{table*}

\section{Proper Orthogonal Decomposition}
The POD is the most widely used technique to reduce data and extract an optimal set of orthonormal bases. In this formulation, basis functions are constructed by collecting temporal snapshots for the variable of interest during the full model solution. Each snapshot vector, $\tilde{U}_x$, holds the values of the local x-velocity at the center of computational cells. The sampled values at the snapshot $s$ are stored in the matrix $\mathbf{\tilde{U}}_{x_s}$ with $\mathcal{N}$ entries ($\mathcal{N}$ being the number of nodes) to construct the snapshot matrices of $\mathbf{\tilde{U}_x}=\left({\tilde{U}}_{x_1}, \ldots, {\tilde{U}}_{x_s}, \ldots, {\tilde{U}}_{x_S}\right)$. For the sake of simplicity, the details will be provided for a general snapshot matrix $\varphi$.

The goal of POD is to find a set of orthogonal basis functions $\left\{\phi_s\right\}, s \in\{1,2, \ldots, S\}$, such that it maximizes:
\begin{equation}
\label{eq_pod1}
\begin{aligned}
& \frac{1}{S} \sum_{s=1}^S\left|<\varphi, \phi_s>_{L^2}\right|^2,
\end{aligned}
\end{equation}
subject to:

\begin{equation}
\label{eq_pod2}
\left\langle\phi_i, \phi_j\right\rangle_{L_2}=\delta_{i j}
\end{equation}
where $\langle\cdot, \cdot\rangle_{L^2} $ is the canonical inner product in L2 norm.

The approach introduced by Sirovich et al.\cite{sirovich1987turbulence} seeks to find an optimal set of basis functions $\phi$ for the optimization problem, i.e. Eq. \eqref{eq_pod1}. This requires performing a Singular Value Decomposition (SVD) of the snapshot matrix $\varphi$ given in the form,

\begin{equation}
\label{eq_pod3}
\begin{aligned}
& \varphi = U \Sigma V^T.
\end{aligned}
\end{equation}

The terms $U \in R^{\mathcal{N} \times \mathcal{N}}$ and $V \in R^{\mathcal{S} \times \mathcal{S}}$ are the matrices that consist of the orthogonal vectors for $\varphi \varphi^{T}$ and $\varphi^{T}\varphi$, respectively and $\Sigma$ is a diagonal matrix of size $\mathcal{N} \times \mathcal{S}$. The non-zero values of $\Sigma$ are the singular values of $\varphi$, and these are assumed to be listed in order of their decreasing magnitude \cite{xiao2015non}. Further on, the POD vectors can be defined as the column vectors of the matrix $U$. These vectors are considered to be optimal in the sense that no other rank $N$ set of basis vectors can be closer to the snapshot matrix $\varphi$ in the Frobenius norm. 

A reduced-order approximation of the field can be described as follows:

\begin{equation}
\label{eq_pod4_5}
\varphi \approx \sum_{n=1}^N \alpha_n(t) \cdot \phi_n.
\end{equation}

The loss of information due to the truncation of the POD expansion set to $N$ vectors can be quantified as follows:

\begin{equation}
\label{eq_pod5}
\begin{aligned}
& L=\frac{\sum_{n=1}^N \lambda_n^2}{\sum_{n=1}^S \lambda_n^2},
\end{aligned}
\end{equation}
where $\lambda$ denotes singular values.

\section{Results and discussion}
\begin{figure*}[htbp]
\nolinenumbers
\nolinenumbers
    \centering
    \begin{subfigure}{1\textwidth}
        \includegraphics[scale=0.075]{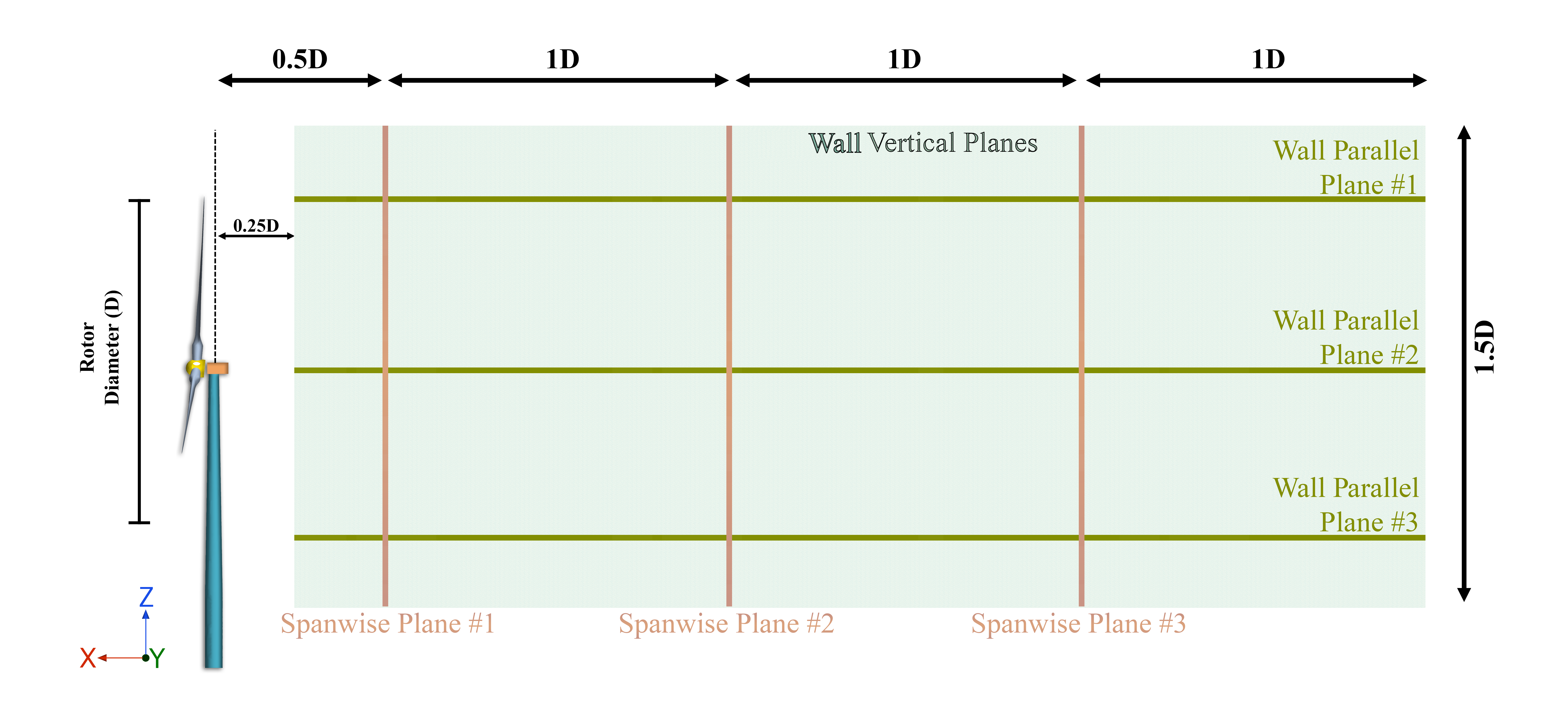} 
        \caption{View normal to the Y-coordinate}
    \end{subfigure}
    
    \begin{subfigure}{1\textwidth}
        \includegraphics[scale=0.075]{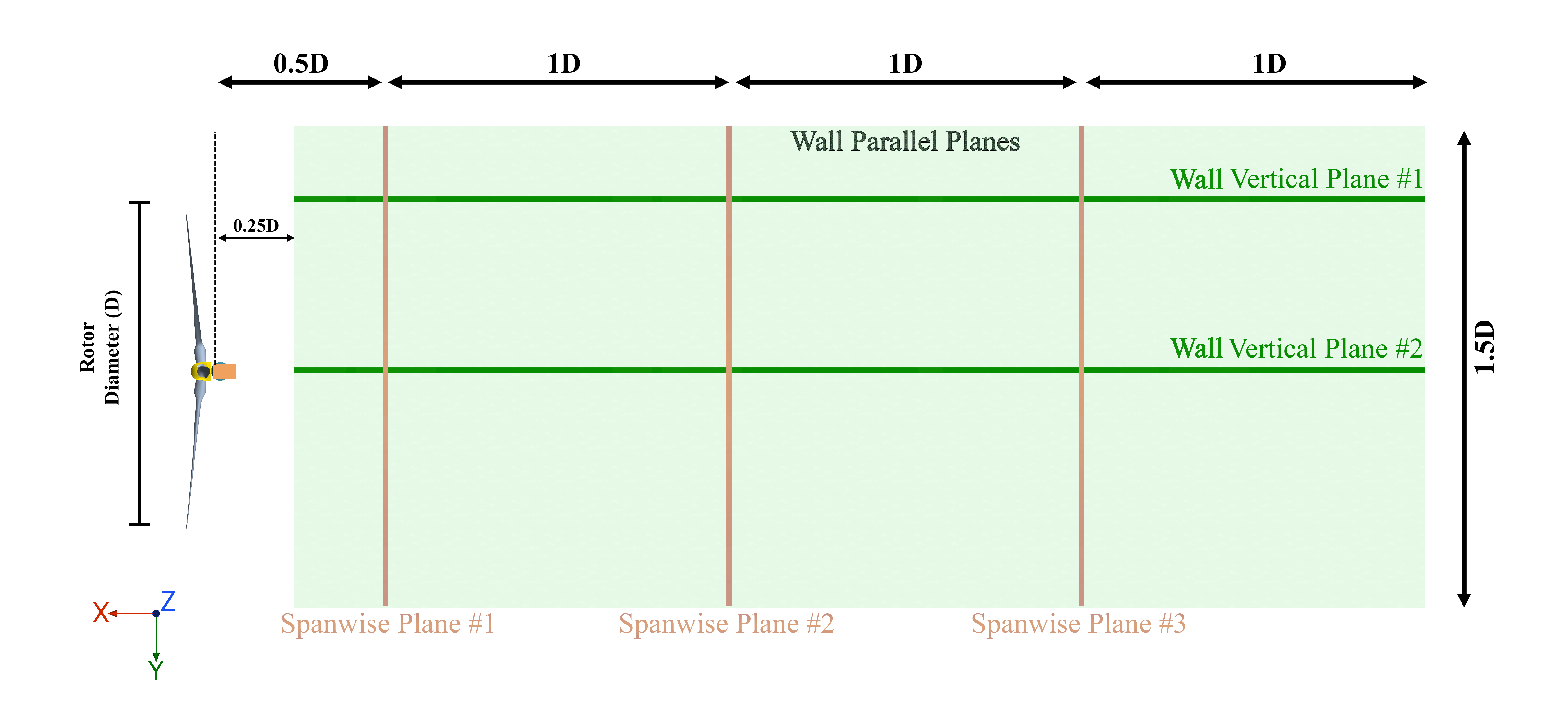}
        \caption{View normal to the Z-coordinate}
    \end{subfigure}

    \caption{\textcolor{darkyellow}{Illustration of the planes for the wake analysis}}
    \label{Fig_ModePlanes}
\end{figure*}
To analyse the wake region behind the turbine a couple of planes were created in the wall-vertical and spanwise directions. As illustrated in Fig. \ref{Fig_ModePlanes}, multiple planes were strategically positioned for detailed wake analysis. In \textcolor{darkyellow}{wall-parallel} (XY) direction, three planes were constructed: one immediately behind the nacelle, another at the turbine tip's highest z-coordinate (+D/2), and a third at the tip's lowest z-coordinate (-D/2) to account for tower effects. Additionally, two planes were created along the \textcolor{darkyellow}{wall-vertical} (ZX) direction: one aligned with the nacelle and another at the farthest point in the Y direction (+D/2) at the tip intersection. This arrangement facilitated a comprehensive examination of the wake flow field in the wall-vertical direction. For evaluating wake characteristics in the spanwise (YZ) direction, three sets of planes were established, ranging from $0.5\si{m}$ to $2.5\si{m}$. The size info of the mentioned planes is presented in Fig. \ref{Fig_ModePlanes}. 

\begin{figure}[htbp]
\nolinenumbers
	\centering
	\includegraphics[width=\linewidth]{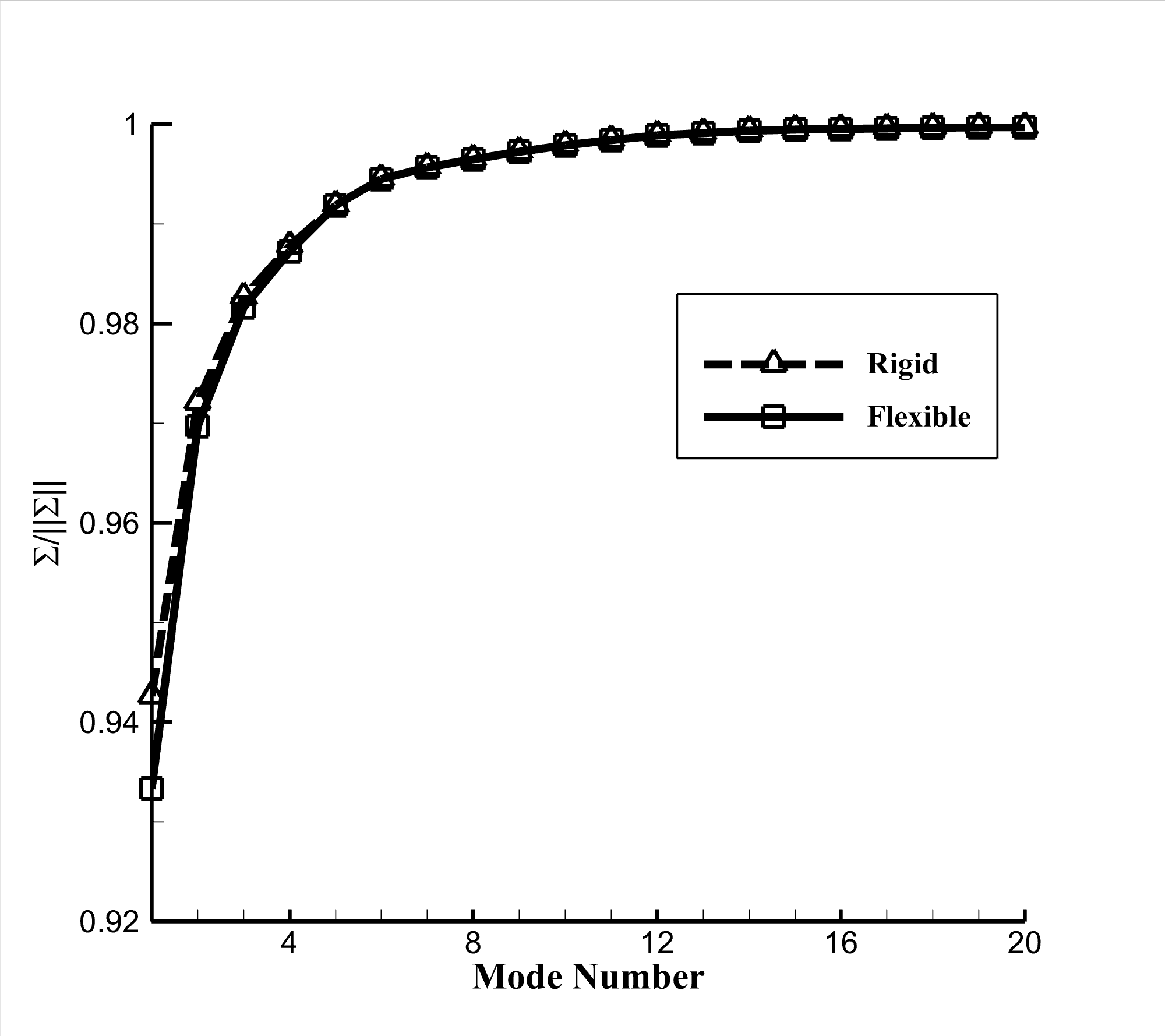}
	\caption{Cumulative modal energy of the rigid and flexible turbines at the \textcolor{darkyellow}{wall-vertical} (ZX) plane behind the nacelle}
	\label{Fig_Modes_XY}
 \end{figure}

\begin{figure}[htbp]
\nolinenumbers
	\centering
	\includegraphics[width=\linewidth]{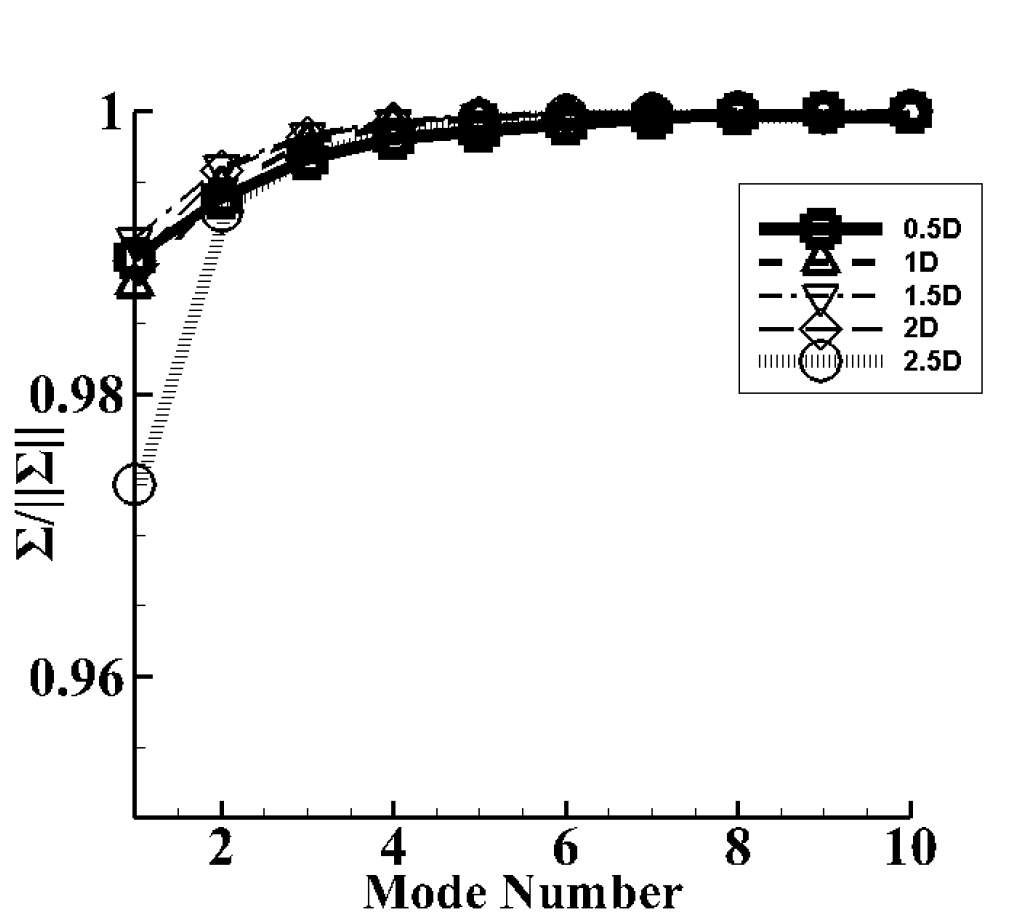}
	\caption{Cumulative modal energy curve at different spanwise (YZ) planes behind the turbine}
	\label{Fig_Modes_Accu}
\end{figure}

\begin{figure*}[t]
\nolinenumbers
	\centering
	\includegraphics[scale=0.3]{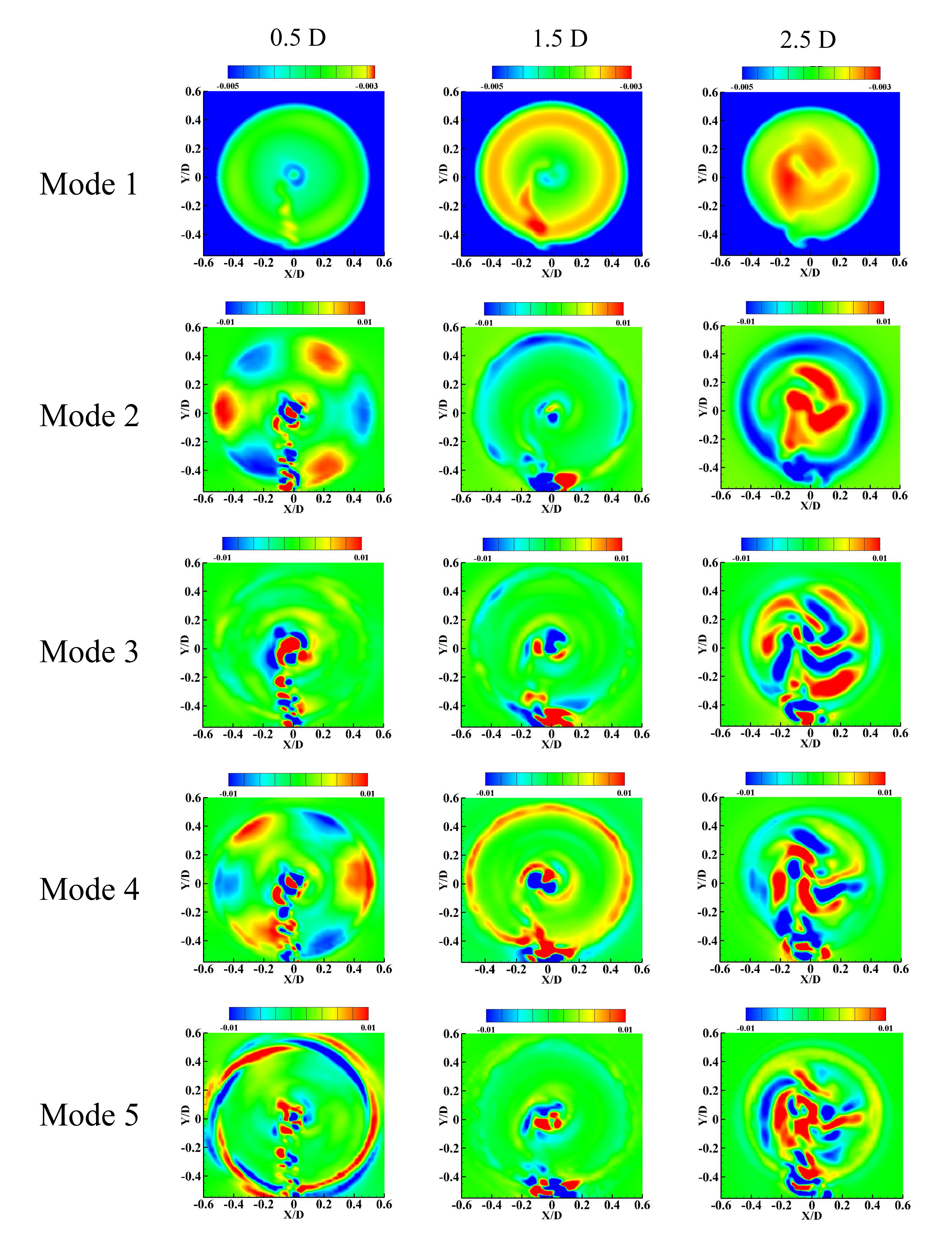}
	\caption{POD modes of velocity at different spanwise (YZ) planes behind the turbine with \textcolor{darkyellow}{flexible blades}}
	\label{Fig_Modes_ZX}
\end{figure*}

\begin{figure*}[t]
\nolinenumbers
	\centering
	\includegraphics[scale=.8]{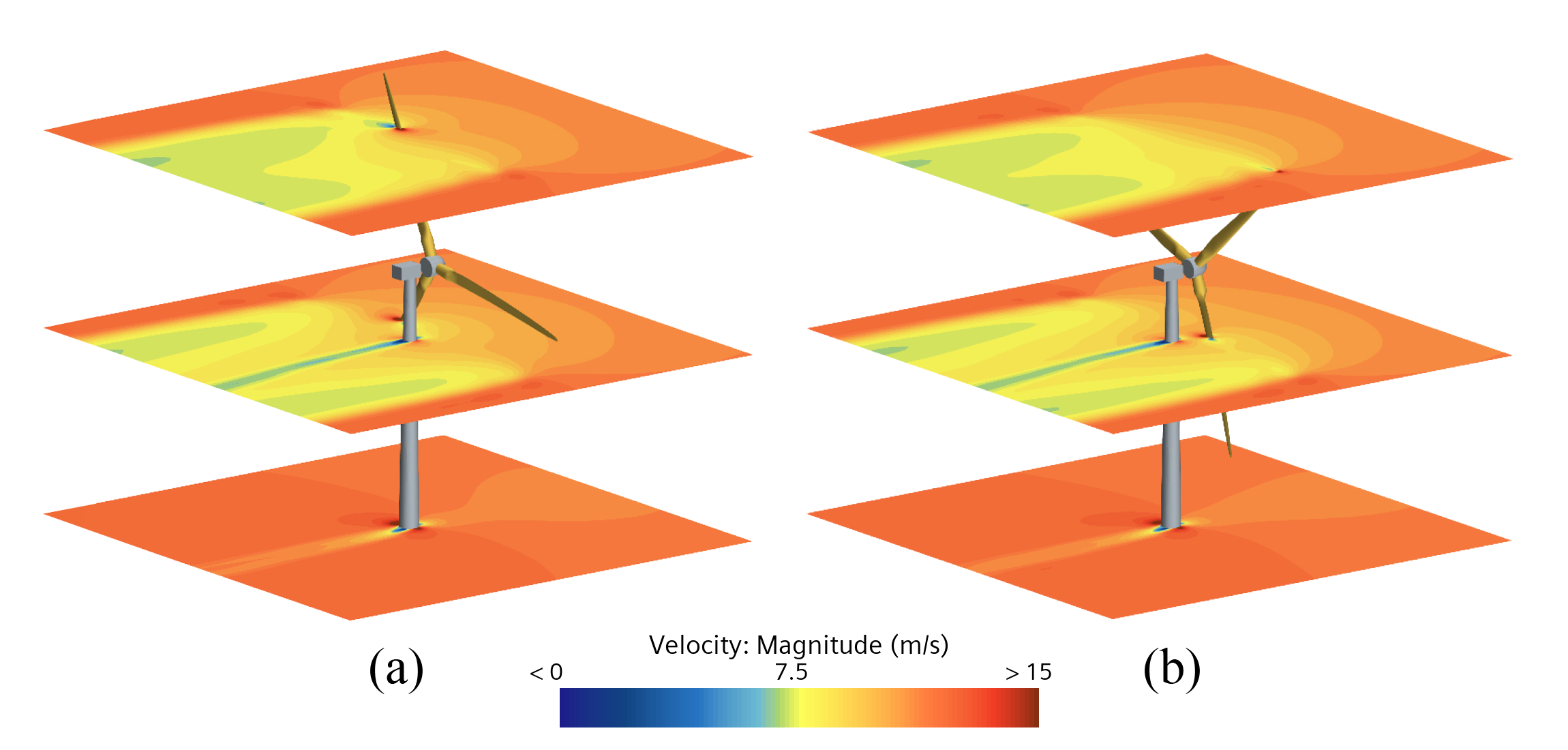}
	\caption{\textcolor{darkyellow}{Contours of velocity magnitude at different wall-parallel (XY) planes for the flexible turbine}}
	\label{Fig_Velo}
\end{figure*}

The modal accumulative energy of the velocity magnitude at the wall-vertical plane is shown in Fig. \ref{Fig_Modes_XY}. Wake regions of rigid and flexible blades demonstrated a similar trend and close values, especially after the second mode. The highest difference was observed within the first mode which corresponds to less than 2\%. It is important to highlight that the first mode for both cases demonstrates a higher difference with the consecutive ones compared to the other modes. More importantly, the first mode of the flexible case has a lower value than that of the rigid indicating that more vortical structures exist in the former. However, the small difference in the first modes of both cases insinuates that difference cannot be fundamentally different and effective in changing the overall behaviour of the wake region. Moreover, the value of the first mode of the flexible turbine was less than that of the rigid turbine due to the blade deflection (harmonic nature of flap-wise and torsional deflections) and corresponding flow alterations. After all, this result indicates the unnecessity of considering blade deformation for wake analysis behind large-scale wind turbines. 

%

The modal cumulative energy at different spanwise planes from 0.5D to 2.5D away from the rotor has been presented in Fig. \ref{Fig_Modes_Accu}. First of all, it should be noted that as the first mode contains the highest energy its color bar scale was chosen to be of a different of range than that of the other modes as also performed in other studies \citep{premaratne2022proper, cherubini2021data, raibaudo2022pod}. It is evident that the farthest plane located 2.5D away from the rotor contains the lowest amount of cumulative energy in the first mode. As a result, it is expected to observe more flow structures at the furthest plane. In addition, the modal energy for the planes from 0.5D to 2D demonstrated a similar trend to each other. Therefore, the flow structures at the near plates could represent analogous vortical characteristics. 

In this regard, in order to have a deeper insight into the wake, the spatial modes of the wall-vertical component of velocity at planes located 0.5D, 1.5D and 2.5D behind the rotor have been presented in Fig. \ref{Fig_Modes_ZX}. The first modes as seen, mainly show the area affected by the presence of the rotor which was also found in the work of Cillis et al.\citep{de2021pod}. Additionally, this area is subtly affected by the presence of the tower and nacelle. To illustrate, in the first mode at 0.5D presented in Fig. \ref{Fig_Modes_ZX}, the affected area by the energy extraction is distorted towards the lower part due to the wake shedding of the tower whose trace can be more clearly seen in the consecutive planes. For instance, in the first mode at 1.5D, the wake shed from the tower becomes more pronounced indicating that moving away from the rotor, the effect of energy extraction is outweighed by the wake shed. Furthermore, the first mode at 1.5D vividly shows the tip vortices emanating from the blades' loading. The decay of this trend can be observed at farther planes at 2.5D. Moreover, the flow structures due to the wake shed from the tower tend to move towards the center and merge with that of the nacelle leading to an integral vortical structure 2.5D away from the rotor. In addition, the symmetrical shape of the first mode starts to fade away with the increase in the distance from the rotor proving that while the work extraction is the main characteristic of the planes near the rotor, the effects of wake shed and tip vortices rise into importance in the farther planes. 

\begin{figure}[htbp]
\nolinenumbers
	\centering
	\includegraphics[scale=0.7]{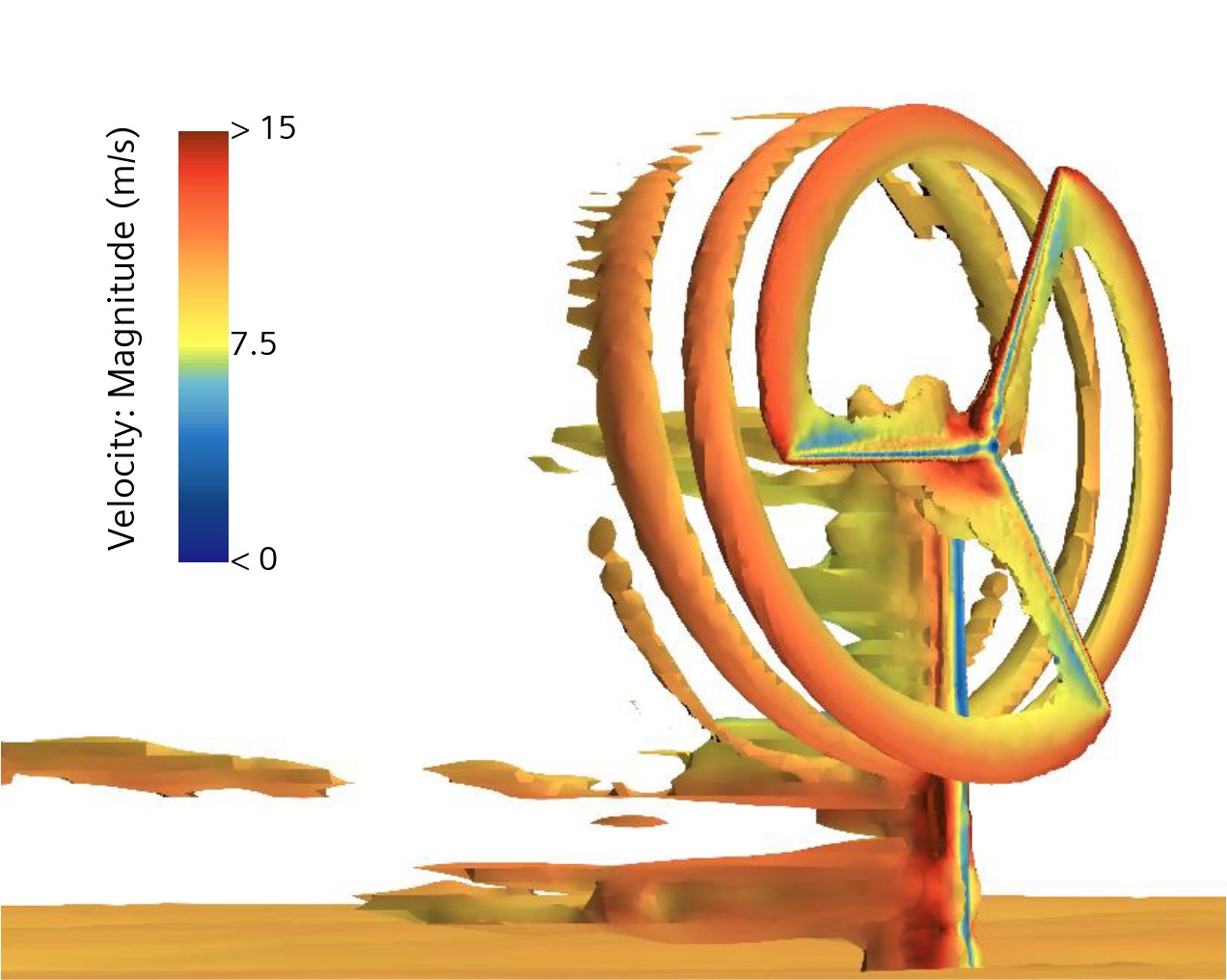}
	\caption{\textcolor{darkyellow}{Iso-surfaces of Q-criterion coloured by velocity magnitude for the flexible turbine}}
	\label{Fig_Vort}
\end{figure}

\begin{figure}[t]
\nolinenumbers
    \centering
    
    \begin{subfigure}{0.45\textwidth}
    \nolinenumbers
        \includegraphics[width=\linewidth]{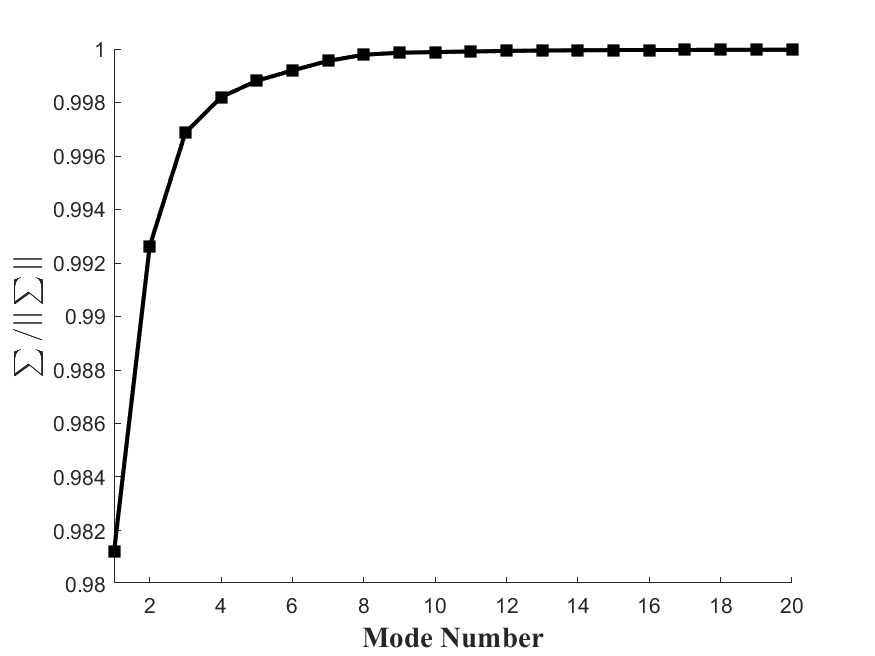}
        \caption{The \textcolor{darkyellow}{wall-vertical} (ZX) plane with $+D/2$ distance from the nacelle}
        \label{acc_stream}
    \end{subfigure}
    
    \begin{subfigure}{0.45\textwidth}
    \nolinenumbers
        \includegraphics[width=\linewidth]{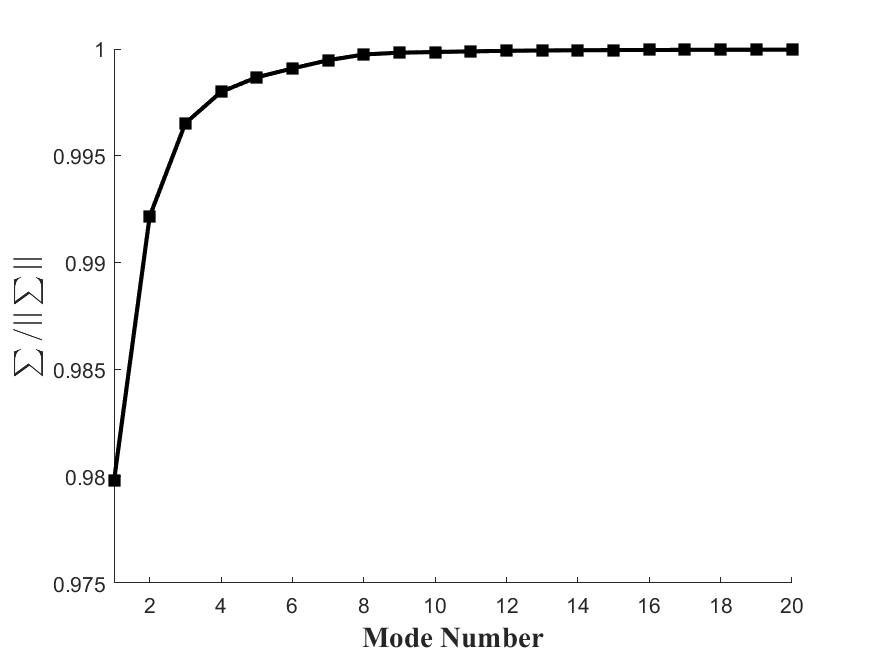} 
        \caption{The \textcolor{darkyellow}{wall-parallel} (XY) plane with $+D/2$ distance from the nacelle}
        \label{acc_span}
    \end{subfigure}

    \caption{Accumulative modal energy on (a) \textcolor{darkyellow}{wall-vertical and (b) wall-parallel} planes \textcolor{darkyellow}{for the flexible turbine}}
    \label{acc_energ_comp}
\end{figure}

\begin{figure}[t]
\nolinenumbers
    \centering
    
    \begin{subfigure}{0.45\textwidth}
    \nolinenumbers
        \includegraphics[width=\linewidth]{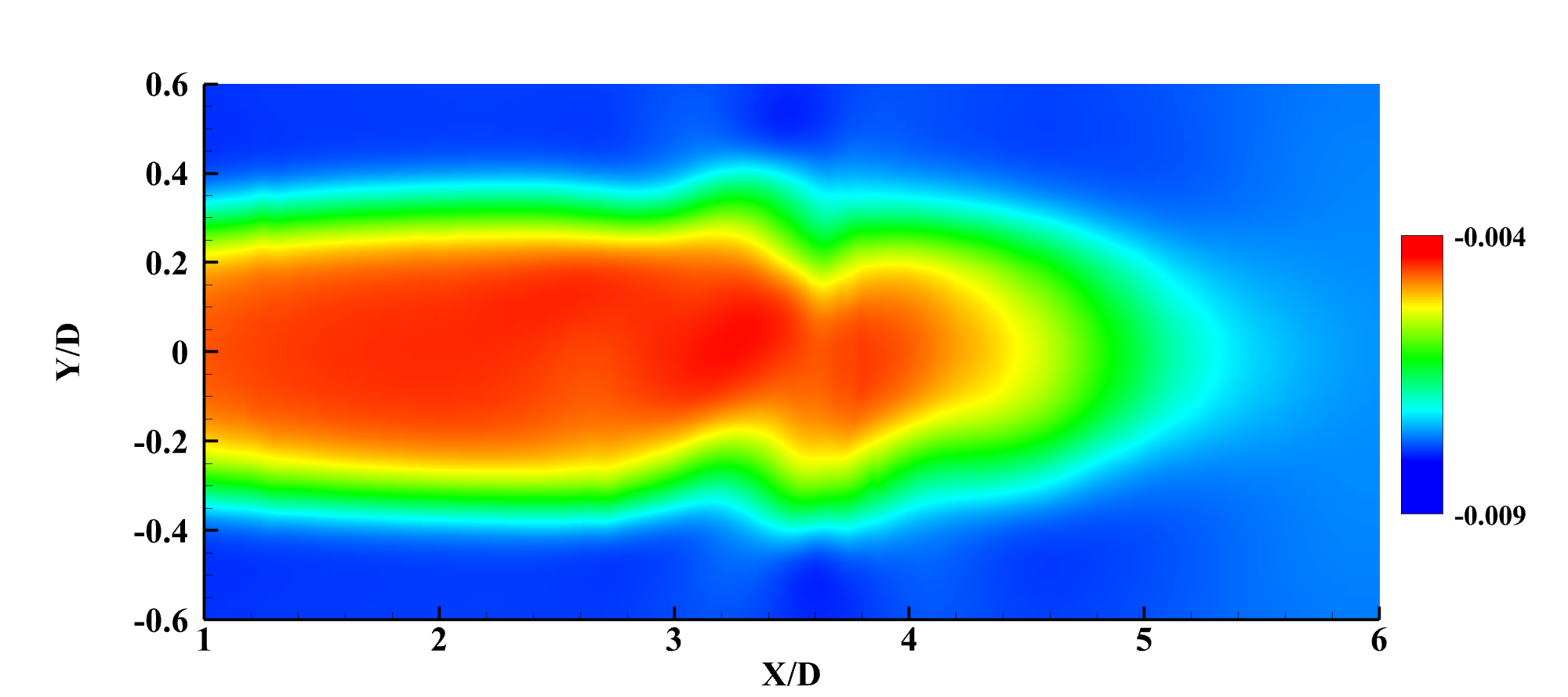}
        \caption{mode 1}
        \label{tip_mode1}
    \end{subfigure}
    
    \begin{subfigure}{0.45\textwidth}
    \nolinenumbers
        \includegraphics[width=\linewidth]{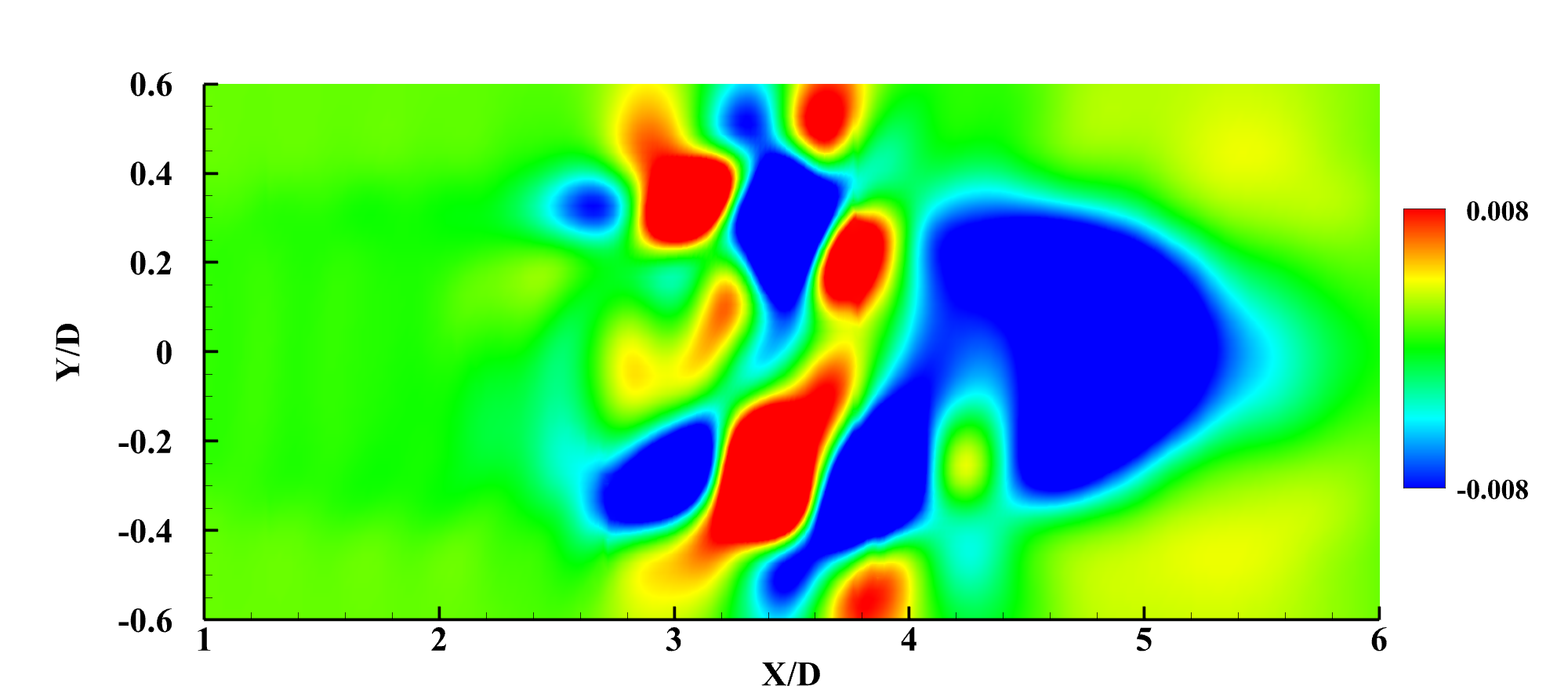} 
        \caption{mode 2}
        \label{tip_mode2}
    \end{subfigure}

    \begin{subfigure}{0.45\textwidth}
    \nolinenumbers
        \includegraphics[width=\linewidth]{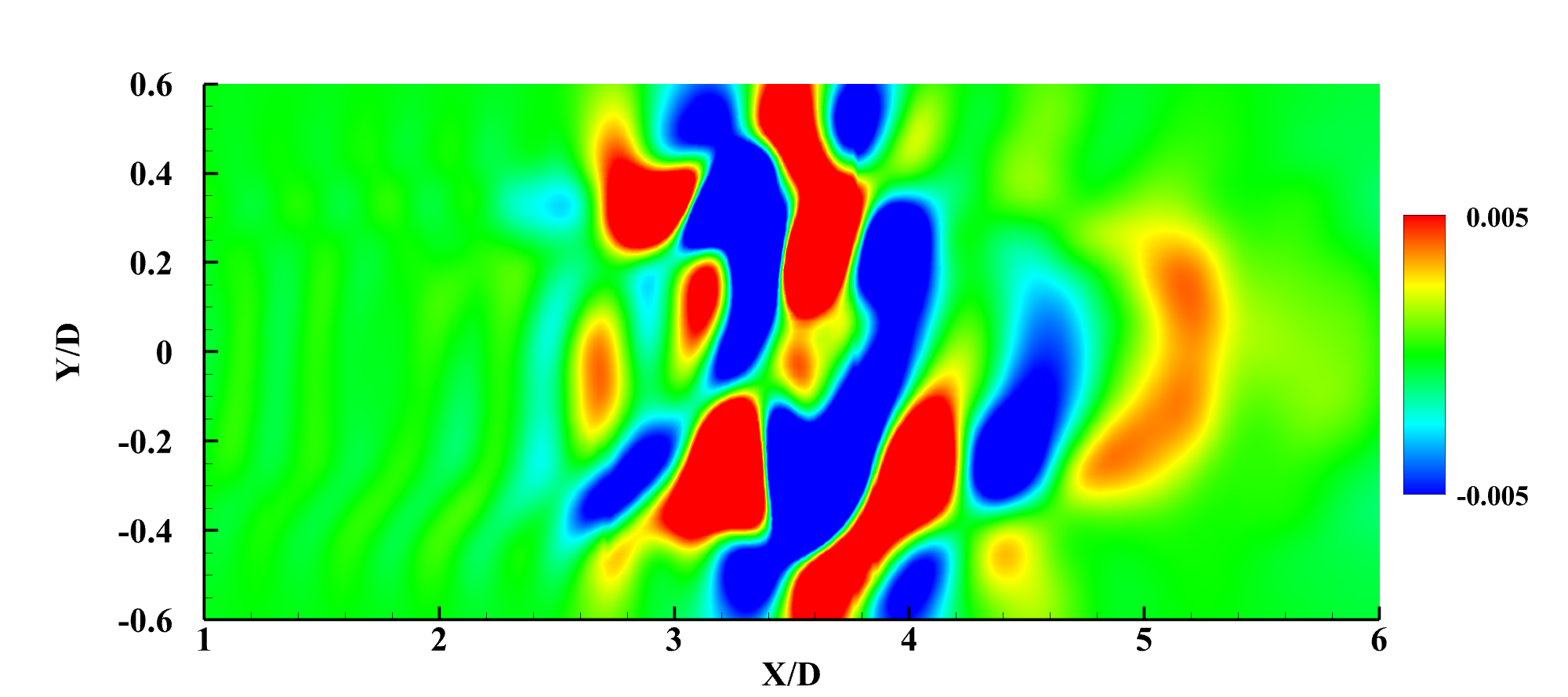} 
        \caption{mode 3}
        \label{tip_mode3}
    \end{subfigure}

    \caption{POD modes on the \textcolor{darkyellow}{wall-parallel} (XY) plane on the blade tip ($+D/2$) \textcolor{darkyellow}{for the flexible turbine}}
    \label{XY_tip}
\end{figure}

\begin{figure}[t]
\nolinenumbers
    \centering
    
    \begin{subfigure}{0.45\textwidth}
    \nolinenumbers
        \includegraphics[width=\linewidth]{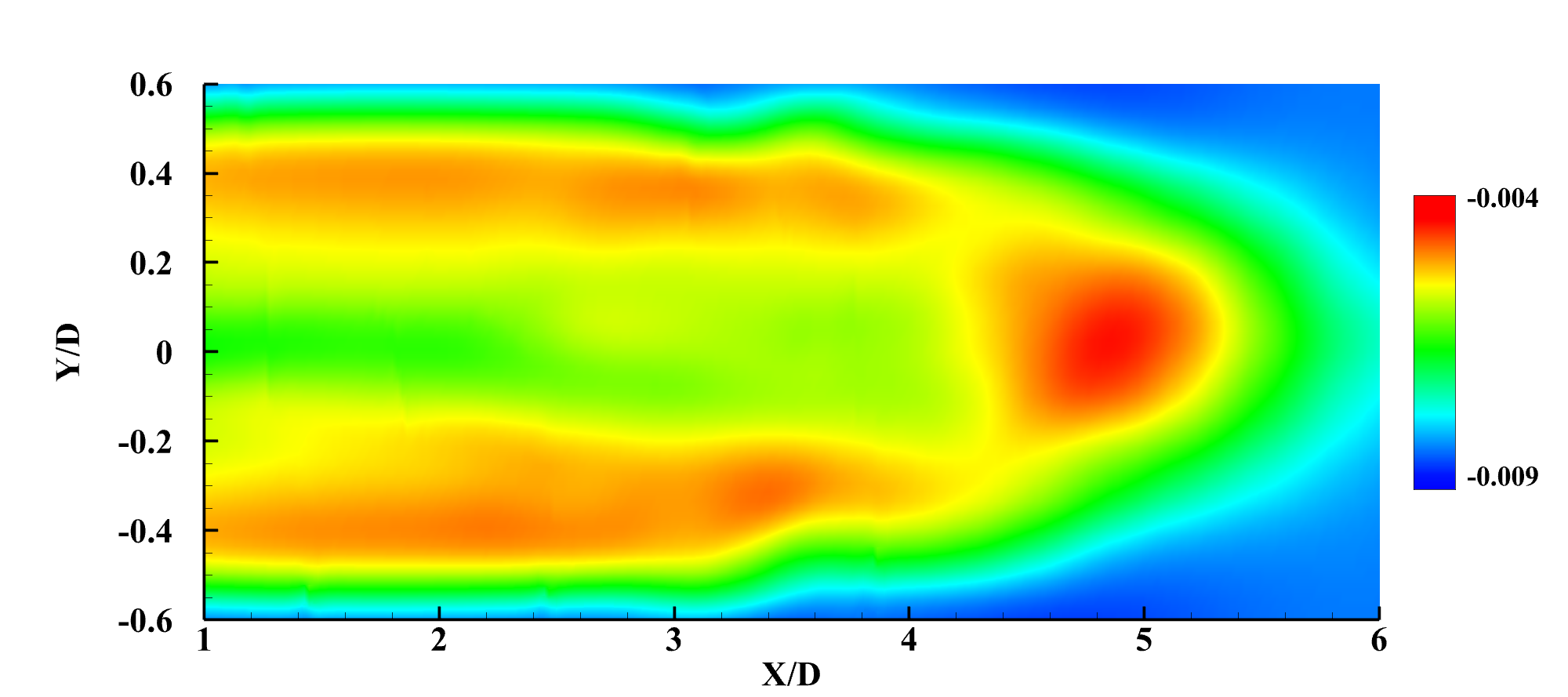}
        \caption{mode 1}
        \label{nacelle_mode1}
    \end{subfigure}
    
    \begin{subfigure}{0.45\textwidth}
    \nolinenumbers
        \includegraphics[width=\linewidth]{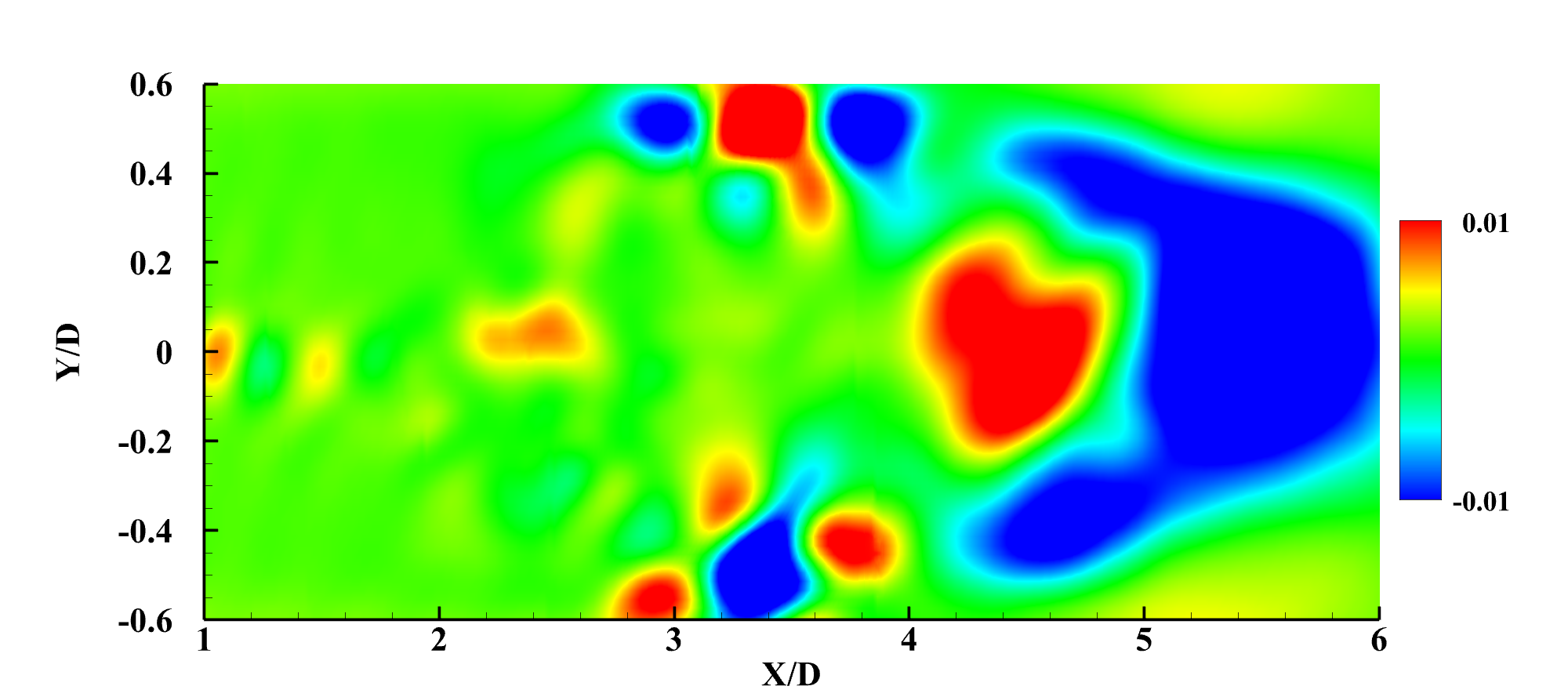} 
        \caption{mode 2}
        \label{nacelle_mode2}
    \end{subfigure}

    \begin{subfigure}{0.45\textwidth}
    \nolinenumbers
        \includegraphics[width=\linewidth]{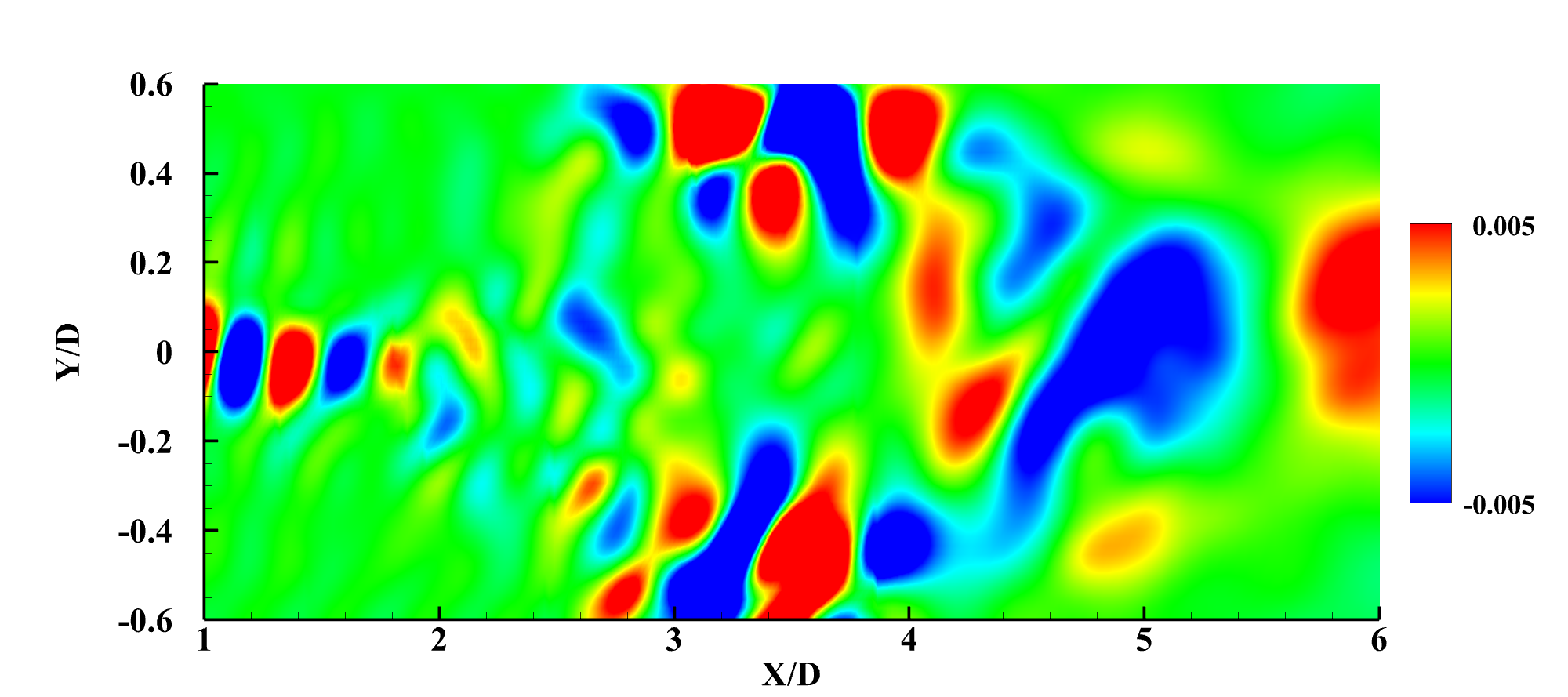} 
        \caption{mode 3}
        \label{nacelle_mode3}
    \end{subfigure}

    \caption{POD modes on the \textcolor{darkyellow}{wall-parallel} (XY) plane behind the nacelle \textcolor{darkyellow}{for the flexible turbine}}
    \label{XY_nacelle}
\end{figure}

\begin{figure}[t]
\nolinenumbers
    \centering
    
    \begin{subfigure}{0.45\textwidth}
    \nolinenumbers
        \includegraphics[width=\linewidth]{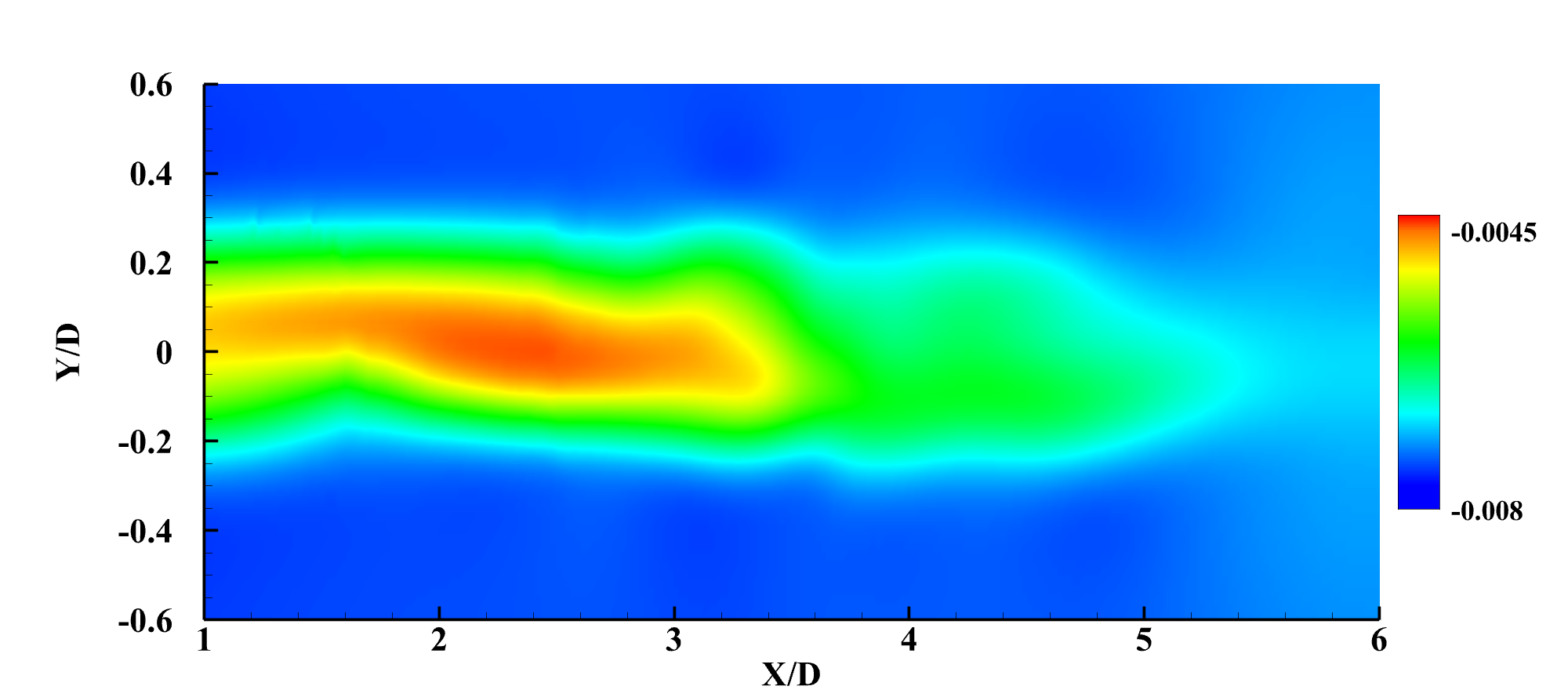}
        \caption{mode 1}
        \label{tower_mode1}
    \end{subfigure}
    
    \begin{subfigure}{0.45\textwidth}
    \nolinenumbers
        \includegraphics[width=\linewidth]{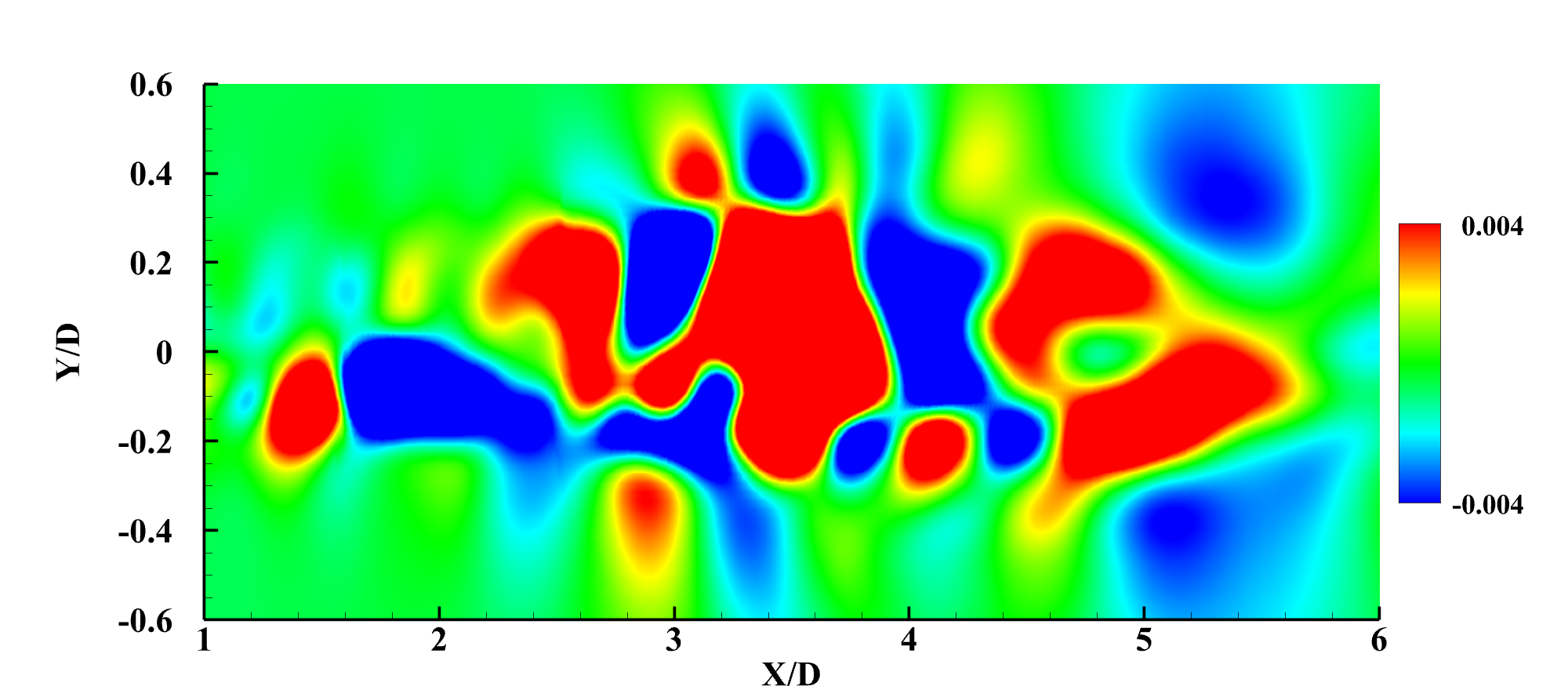} 
        \caption{mode 2}
        \label{tower_mode2}
    \end{subfigure}

    \begin{subfigure}{0.45\textwidth}
    \nolinenumbers
        \includegraphics[width=\linewidth]{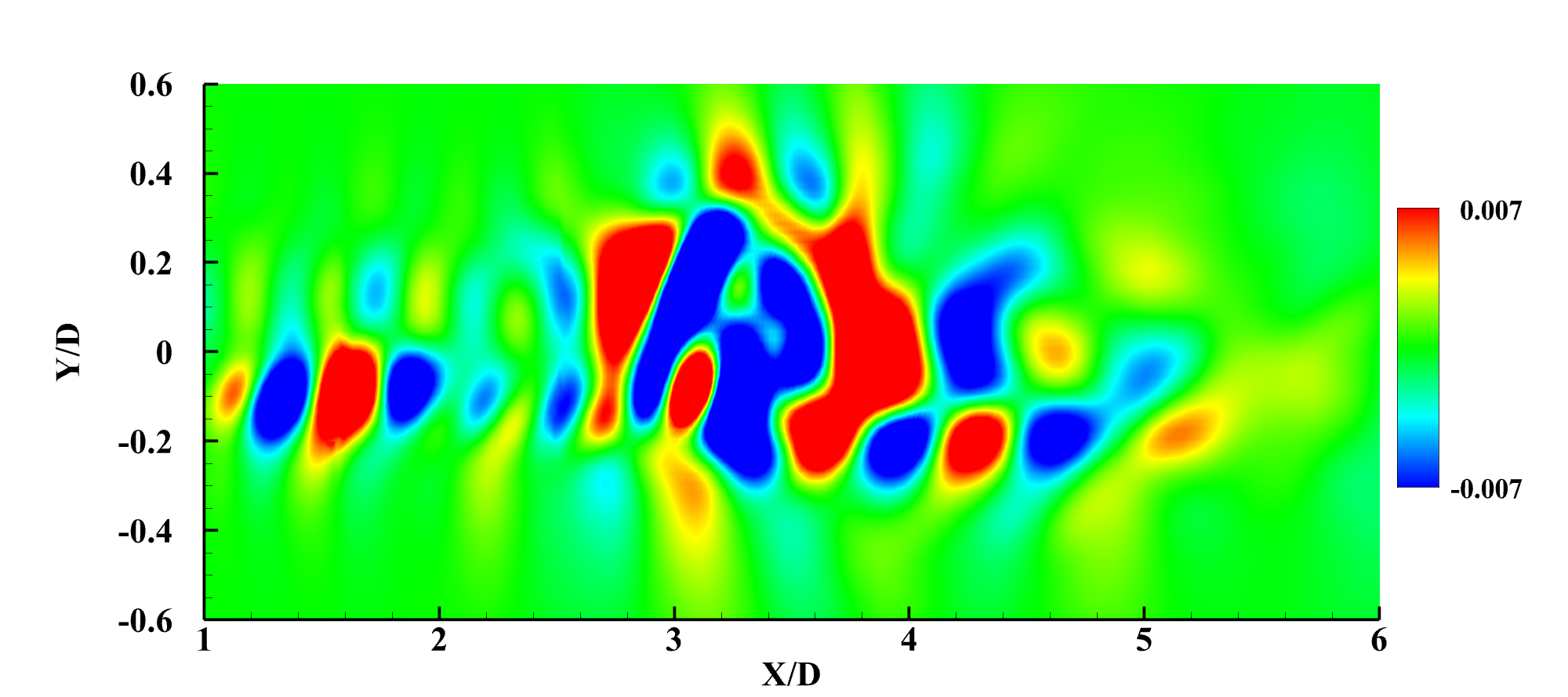} 
        \caption{mode 3}
        \label{tower_mode 3}
    \end{subfigure}

    \caption{POD modes on the \textcolor{darkyellow}{wall-parallel} (XY) plane behind the tower ($-D/2$) \textcolor{darkyellow}{for the flexible turbine}}
    \label{XY_tower}
\end{figure}

\begin{figure}[t]
\nolinenumbers
    \centering
    
    \begin{subfigure}{0.45\textwidth}
    \nolinenumbers
        \includegraphics[width=\linewidth]{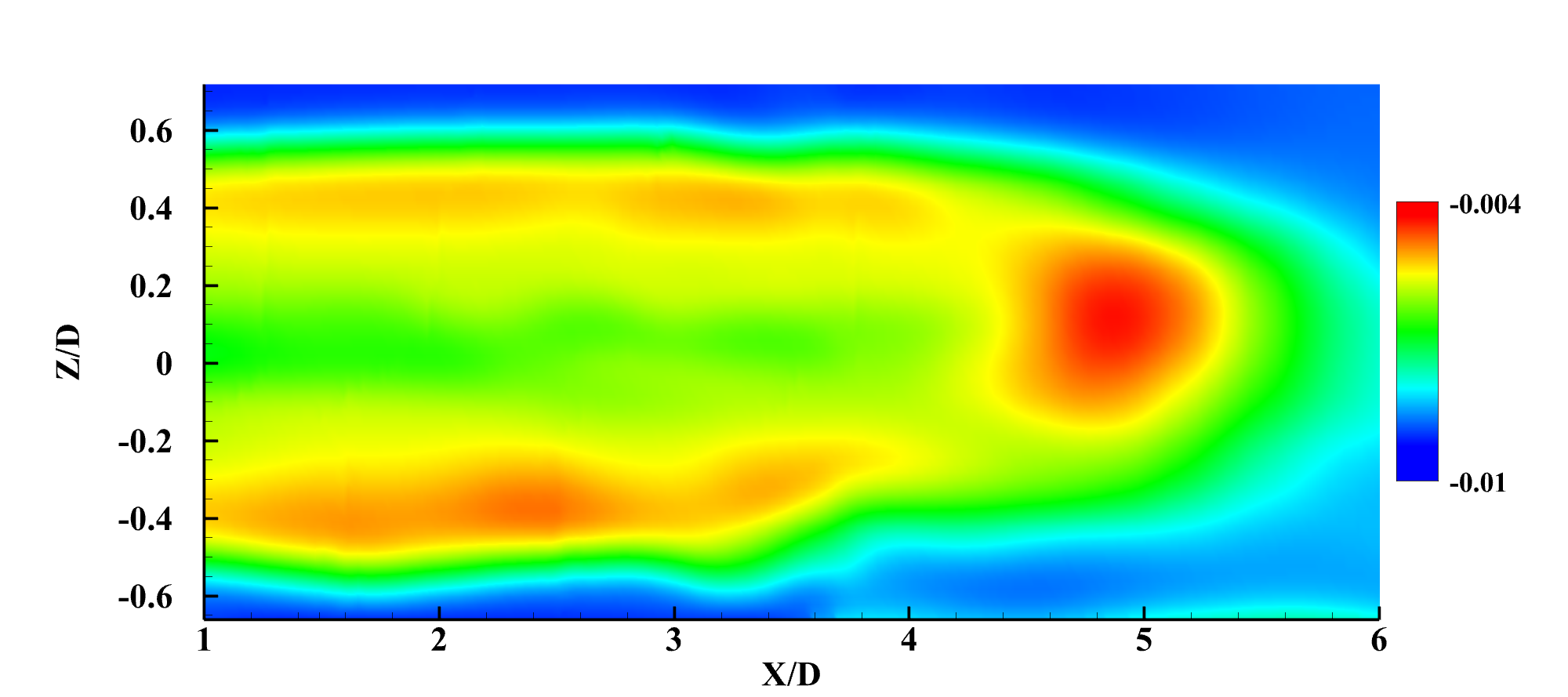}
        \caption{mode 1}
        \label{near_mode1}
    \end{subfigure}
    
    \begin{subfigure}{0.45\textwidth}
    \nolinenumbers
        \includegraphics[width=\linewidth]{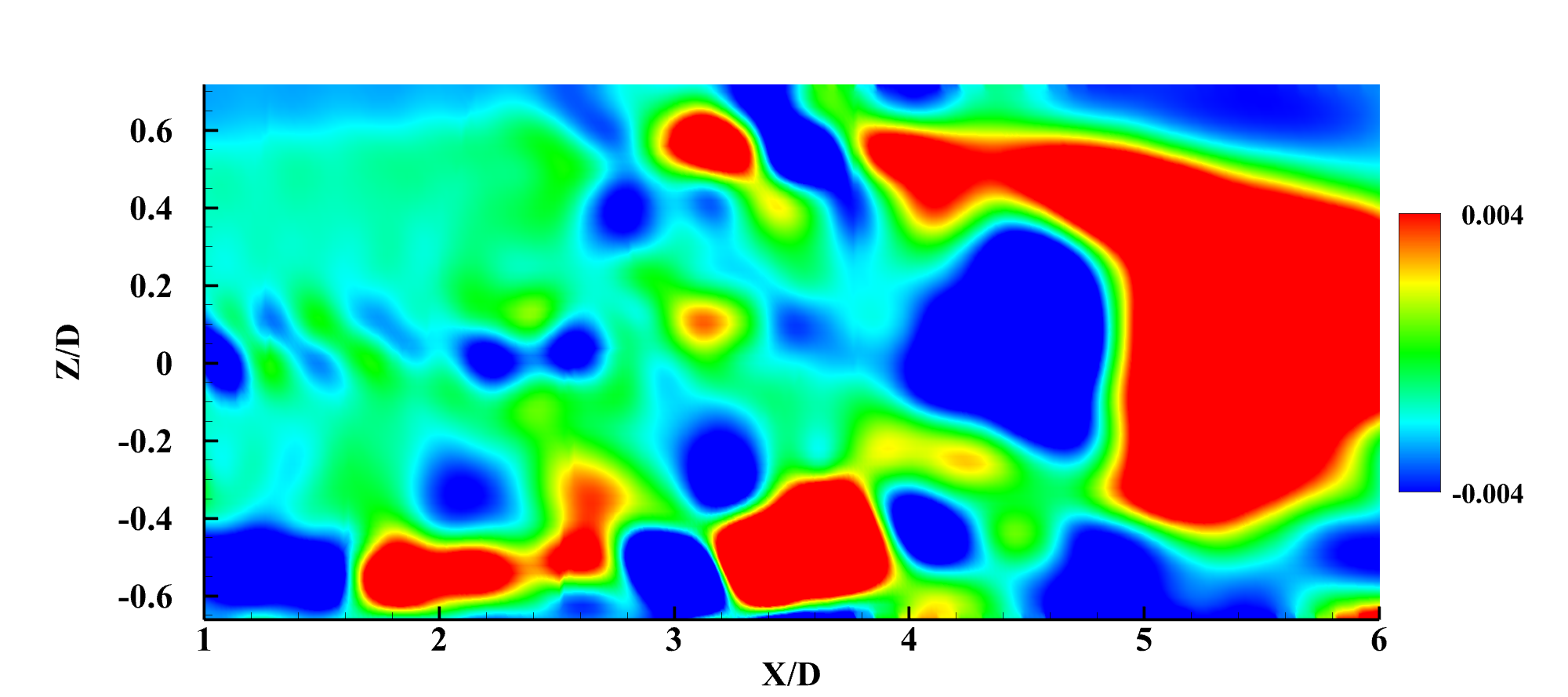} 
        \caption{mode 2}
        \label{near_mode2}
    \end{subfigure}

    \begin{subfigure}{0.45\textwidth}
    \nolinenumbers
        \includegraphics[width=\linewidth]{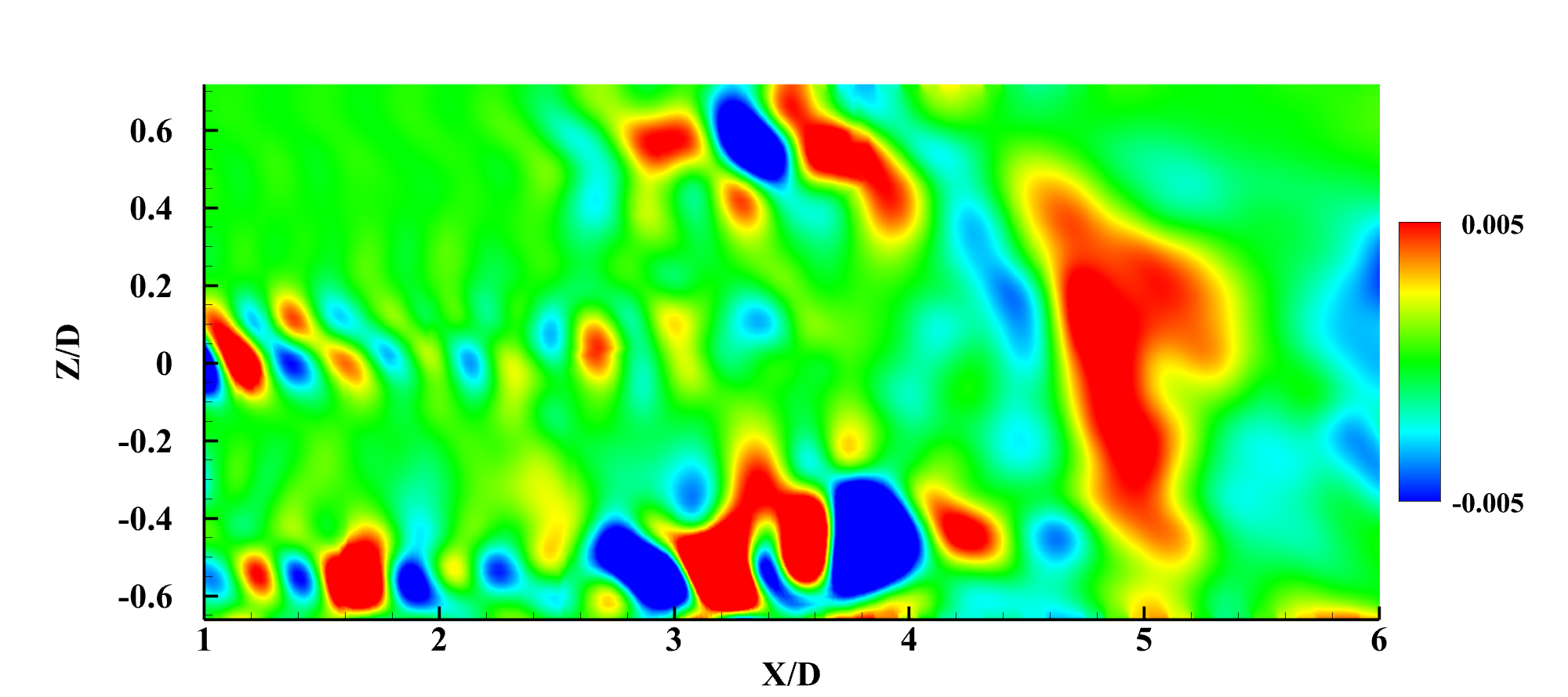} 
        \caption{mode 3}
        \label{near_mode3}
    \end{subfigure}

    \caption{POD modes on the \textcolor{darkyellow}{wall-vertical} (ZX) plane behind nacelle \textcolor{darkyellow}{for the flexible turbine}}
    \label{near_mode}
\end{figure}

\begin{figure}[t]
\nolinenumbers
    \centering
    
    \begin{subfigure}{0.45\textwidth}
    \nolinenumbers
        \includegraphics[width=\linewidth]{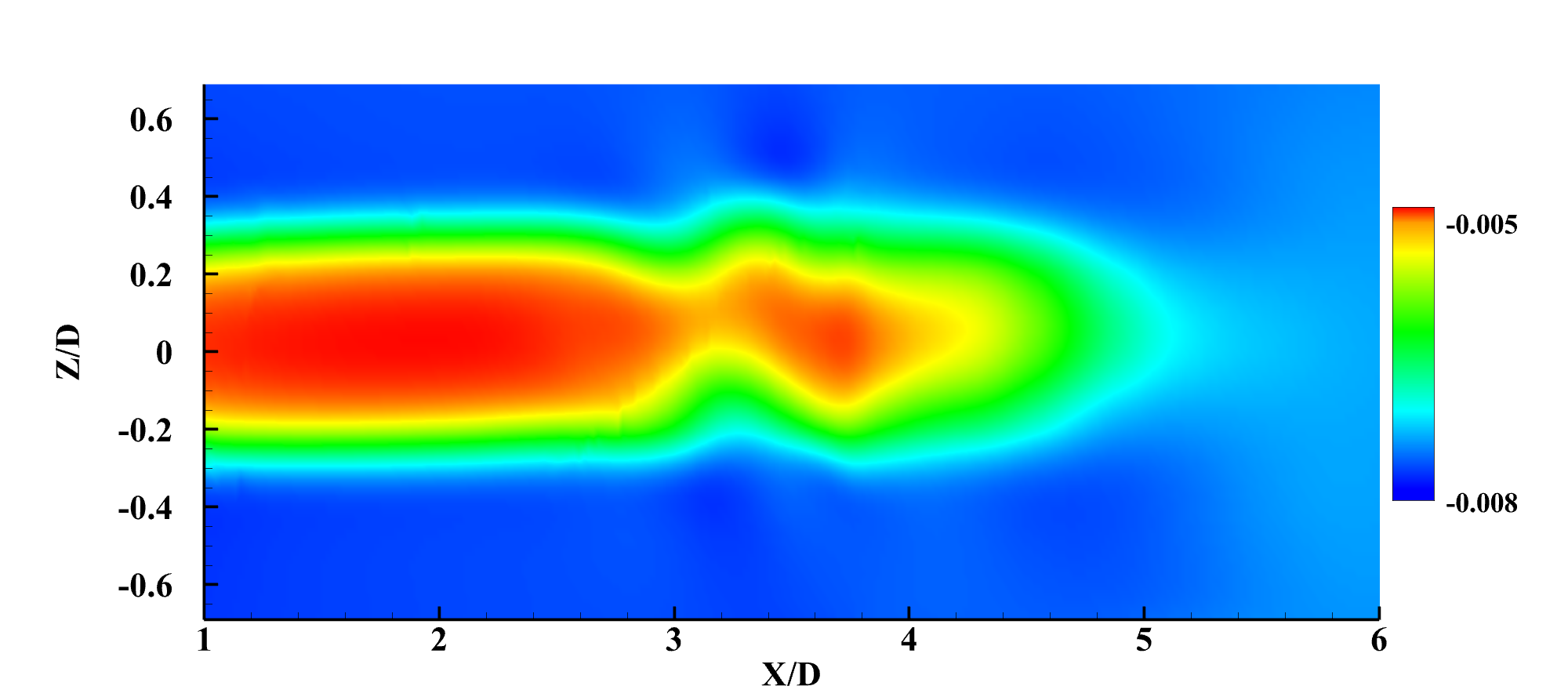}
        \caption{mode 1}
        \label{far_mode1}
    \end{subfigure}
    
    \begin{subfigure}{0.45\textwidth}
    \nolinenumbers
        \includegraphics[width=\linewidth]{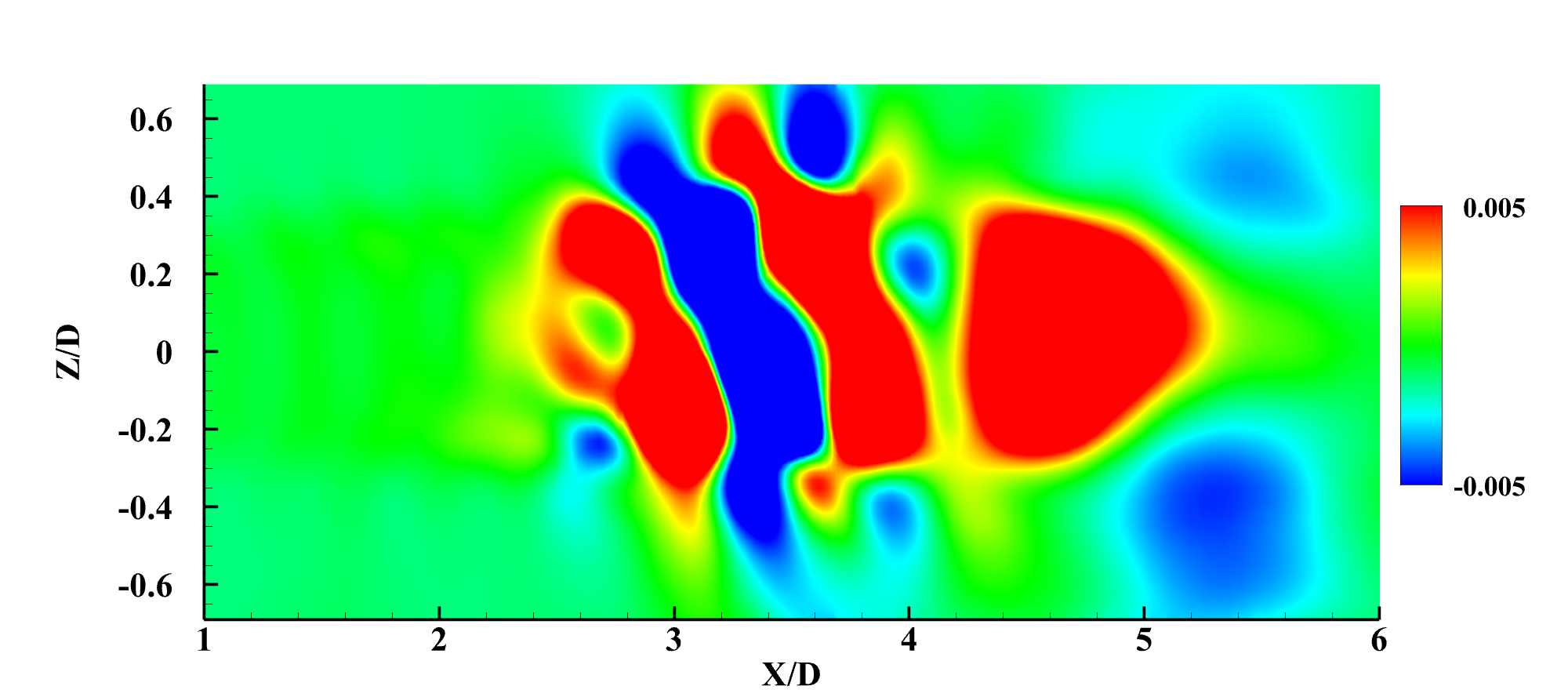} 
        \caption{mode 2}
        \label{far_mode2}
    \end{subfigure}

    \begin{subfigure}{0.45\textwidth}
    \nolinenumbers
        \includegraphics[width=\linewidth]{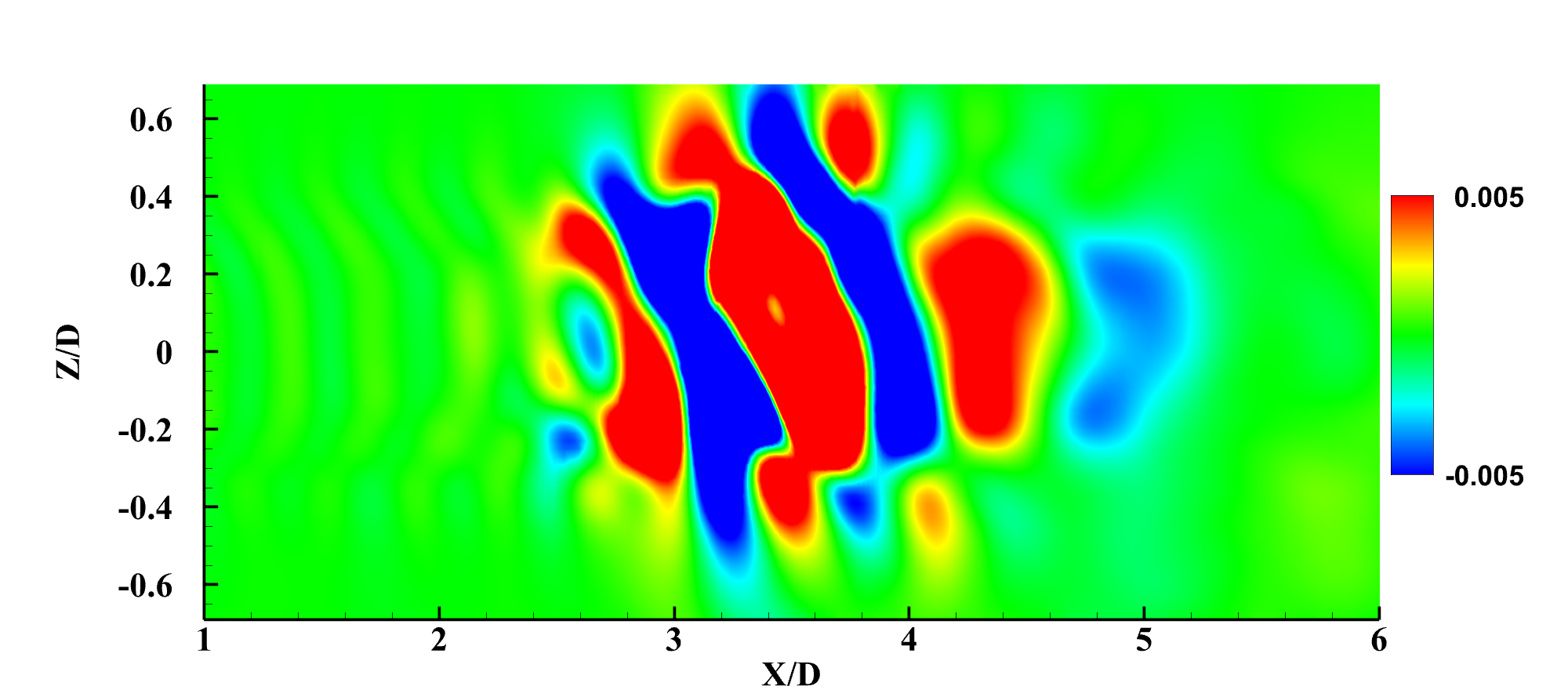} 
        \caption{mode 3}
        \label{far_mode3}
    \end{subfigure}

    \caption{POD modes on the \textcolor{darkyellow}{wall-vertical} (ZX) plane at the blade's tip ($+D/2$) \textcolor{darkyellow}{for the flexible turbine}}
    \label{far_mode}
\end{figure}


The second mode at 0.5D from the rotor illustrates the structures stemming from the tower wake and the tip vortices which are the most effective patterns in the wake of the wind turbine. This observation has also been reported in other studies\citep{andersen2014reduced, iungo2015data}. It should be noted that with the increase in the distance from the wind turbine, the decay of the structures at the plane in 0.5D is observed at the same mode number. To illustrate, the flow structures shedding behind the tower are observed at the same positions, however, at smaller dimensions which are also observed for the tip vortices. Moreover, as also shown for the first modes, the flow structures due to the flow interaction with the tower begin to move towards the central part forming an integral structure at 2.5D behind the rotor. However, it should be noted that the central integral structure at 2.5D at the first and the second modes decays into separate structures in the next modes. In addition, the structures emanating from the flow interaction with the tower prove to sustain themselves for a longer distance and through higher mode numbers. This observation insinuates that the tower wake along with those stemming from the blades plays a more important role in the far wake compared to the tip vortices. Also, it is concluded that the exclusion of the tower leads to missing a major part of the wake structures, especially at far wake positions which have also been reported in the work of Premaratne et al.\citep{premaratne2022proper}, and Wang et al. \citep{wang2012study}.

In order to investigate the effect of the tower on the downstream wake, the velocity magnitude contours at different sections have been presented in Fig. \ref{Fig_Velo}. \textcolor{darkyellow}{It is important to highlight that, given the periodic nature of wake behaviour, the instances selected for analysis in this section were chosen to effectively capture and illustrate the influence of the tower on the wake region. These instances were carefully chosen to best represent and elucidate this behaviour, ensuring a comprehensive investigation into the impact of the tower on wake dynamics. As seen, the wake behind the tower leaves a significantly more sustainable trace compared to the wake of the blades. This is due to the fact that as also evident in the POD modes, the blades at the rated TSR undergo a suitable angle of attack leading to a less intense wake shed compared to the tower.} As a result, it is found that the main elements affecting the wake of the wind turbine are the tower, tip vortex and to a lesser extent the blades especially, at the mid and tip sections. As for the comparison of the tower wake and the blade tip vortices which was also drawn in the spatial modes section, the vorticity field coloured by the velocity magnitude has been presented in Fig. \ref{Fig_Vort}. It is observed that the tip vortices sustain themselves for a longer distance compared to the blade wake, while as also previously explicated the tower causes the most significant effect on the wake region which could be traced for a much longer distance. In addition, the merging of the blade wake with that of the tower in the relatively farther planes behind the rotor could be observed which was explicitly discussed in the spatial mode section.   

In order to gain a deeper understanding of the wake, the accumulative energy and the POD modes of the vertical component of velocity in the wall-vertical and wall-parallel planes both on the blade's tip have been presented in Fig. \ref{acc_energ_comp}. First of all, it is evident in the figure that the accumulative energy of both planes demonstrates quite similar trends and values proving that with three modes almost all the dominant flow structures could be extracted. In this regard, the wake behaviour has been deeply investigated on the wall-vertical planes at the blade tip, nacelle and the tower depicted in  Fig.\ref{XY_tip}, Fig.\ref{XY_nacelle} and Fig.\ref{XY_tower}, respectively. 

As seen, the first mode in all the planes corresponds to the area affected by the presence of the body responsible for the energy extraction as also reported in the work of Sheidani et al \cite{sheidani2023assessment}. However, as shown in Fig.\ref{tip_mode1}, Fig.\ref{nacelle_mode1} and Fig.\ref{tower_mode1} that near the blade tip, the wake region is much more uniform compared to the other two planes located at a lower altitude for the first mode. This observation proves that despite the highest deflection which is experienced at the tip in flexible blades, the adversarial effect of the nacelle and tower is yet significantly more noticeable. In addition, the first mode for the plane at the blade tip and the tower bear a closer resemblance with respect to that of the nacelle. As a result, the nacelle is expected to leave more remarkable traces in the wake compared to other parts of the wind turbine. Also, the wake region on the blade tip plane appears to be far more uniform than the planes corresponding to cross sections of lower aerodynamic efficiency i.e. nacelle and tower. 

This finding could be better understood by moving on to the next modes. For instance, the trend in the second modes for all the cross sections in spanwise planes shown in Fig.\ref{tip_mode2}, Fig.\ref{nacelle_mode2} and Fig.\ref{tower_mode2} reveals that the dominant structures correspond to the entrainment stemming from the wake region collapse \cite{sheidani2023assessment2}. However, considering the cross sections of lower aerodynamic efficiencies i.e. at the tower and nacelle, other structures emanating from the flow separation in congruous with sources of performance loss are observed. To illustrate, the second mode at the tower plane demonstrates typical structures in the wake of a cylinder \citep{noack2016recursive} insinuating the importance of the tower on the wake and its modal behaviour. Furthermore, regarding the second mode on the nacelle plane, the structures related to the vortex shedding from the nacelle geometry are analogous to that of the tower plane, nonetheless, at a lower intensity. In this regard, the flow structures of the second mode at the blade tip plane are mainly concentrated around the wake region collapse indicating the higher performance of the tip section as the structures due to the flow separation being of the same modal energy of the expected dominant structure i.e. wake region collapse is not observed. Therefore, considering the effect of flexibility of the blades does not cause the hindrance against the power production by the blade to overweight that of the nacelle and tower. 

Presented Fig.\ref{tip_mode3}, Fig.\ref{nacelle_mode3} and Fig.\ref{tower_mode 3} demonstrate the third POD modes for the blade tip, nacelle and tower sections respectively. The decay of the structures observed previously in the second mode could be observed in all three planes. The third mode for the blade tip and the nacelle show similar symmetrical structures near the wake region borders which are more pronounced near the tip. These structures come into existence due to the velocity difference between the wake region and the area beyond it as explained in \citep{sheidani2023assessment, sheidani2023assessment2}. However, the third mode of the tower plane, as seen in Fig.\ref{tower_mode 3}, shows the structures to be mostly concentrated near the central part of the wake region. Therefore, as expected, the section with the lowest efficiency cannot bring about a considerable velocity gradient between the wake and the far-field regions. It should be noted that, as for the nacelle and tower planes modes, the decay of the structures emanating from the nacelle is still observed, while the flow structures at the blade tip are more pronounced near the wake region collapse proving the higher efficiency of this section.

Fig.\ref{near_mode} and Fig.\ref{far_mode} display the spatial POD modes on the wall-vertical plane, illustrating the potential interactions and phenomena on the wake region. It is important to highlight that two distinct planes are used for evaluation; one intersecting with the nacelle, called the \textit{near plane}, and the other positioned farther from the nacelle at $+D/2$, referred to as the \textit{far plane}. This allows for separate assessments of their effects on the wake. 

Beginning with the first mode on the near plane presented in Fig.\ref{near_mode1}, it is evident that the modal behaviour is quite similar to that of the spanwise plane at the nacelle rather than the tower plane. However, the first POD mode at the far plane exhibits a greater resemblance with the tower plane. This proves that the effect of the nacelle on the wake could be regarded as even more influential than the tower. The significant discrepancies between the modal behaviour of the spanwise and the wall-parallel planes begin after the second mode shown in g.\ref{near_mode2} and Fig.\ref{far_mode2}. To illustrate, the second mode of both near and far planes, despite the presence of similar structures especially, near the point of wake region collapse, the wall-parallel plane suggests more structures emanating mostly from the tower. Therefore, it is found that the effects of the tower on the wake region could be taken into consideration after the second mode in the near plane, specifically. 

The third mode, while sustaining the main trend, shows the decay of the structure observed in the previous mode as expected \cite{sheidani2023assessment}. Nonetheless, the far and near planes reveal a fundamental difference regarding the behaviour of the third mode. As seen in Fig.\ref{far_mode3}, the main structures on the far plane are mainly concentrated around four to five times the wind turbine diameter illustrating the effect of wake region collapse. On the other hand, the third mode on the near plane demonstrates other structures than those corresponding to the decay of the second mode near the wake region. First of all, the vortex shedding of the nacelle around the centre line is evident in Fig.\ref{near_mode3} which is obviously absent on the far plane due to its distant position from the nacelle. 

Moreover, the symmetrical structures near the boundaries of the wake region are observed on the near which is yet another remarkable difference between the far and near plane's third modes. It is observed in Fig.\ref{near_mode3} that the presence of the symmetrical structures is more pronounced on the lower part of the part of the wake region. This is attributed to the fact that the lower part is prone to the flow structures ensuing from the presence of the tower, while the upper part of the wake region is not exposed to a significant source of vortex shedding upstream. As a result, the symmetrical structures on the lower part are intensified resulting in an asymmetrical wake region shape on the near plane. 

\begin{figure}[t]
\nolinenumbers
    \centering

    \begin{subfigure}{0.45\textwidth}
    \nolinenumbers
        \includegraphics[width=\linewidth]{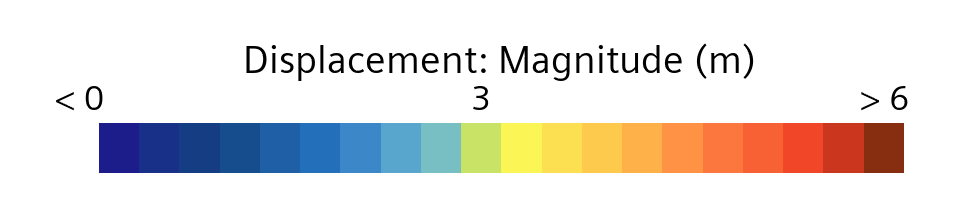} 
    \end{subfigure}
    \begin{subfigure}{0.45\textwidth}
    \nolinenumbers
        \includegraphics[width=\linewidth]{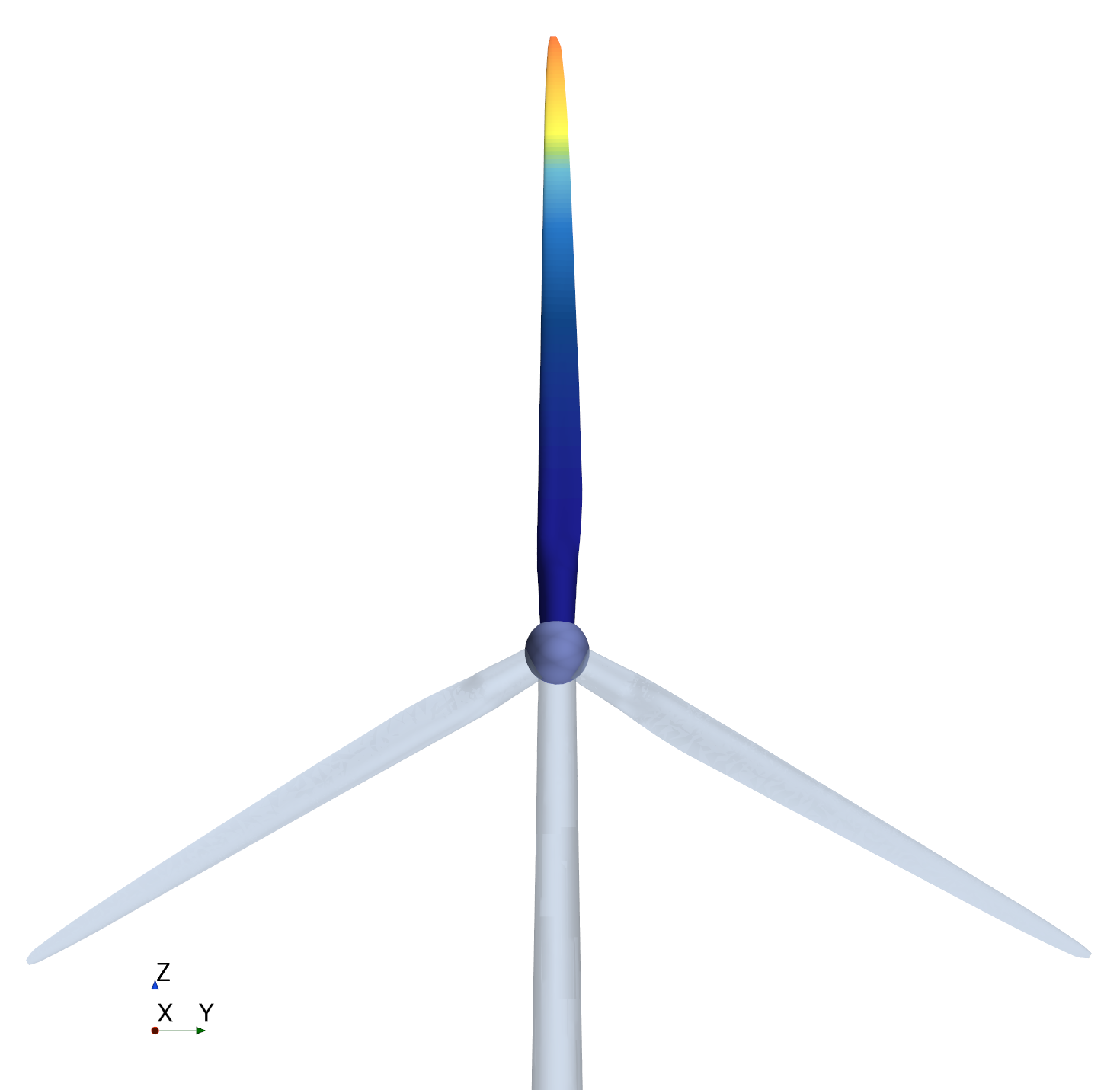} 
        \caption{Highest blade position}
        \label{defl 1}
    \end{subfigure}
    
    \begin{subfigure}{0.45\textwidth}
    \nolinenumbers
        \includegraphics[width=\linewidth]{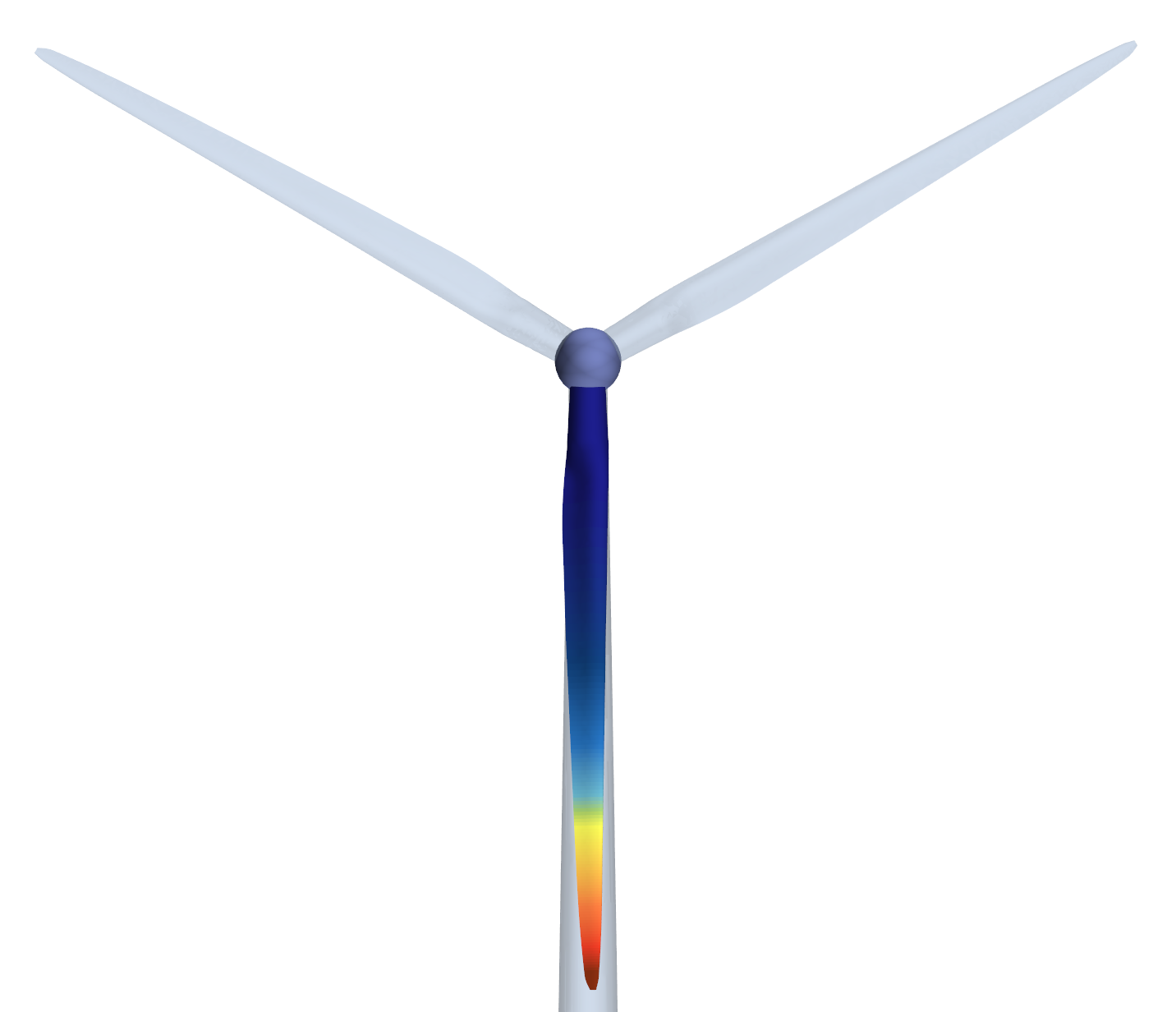}
        \caption{Lowest blade position}
        \label{defl 2}
    \end{subfigure}

    \caption{\textcolor{darkyellow}{Blade deflection at two critical positions}}
    \label{blade_deflection}
\end{figure}

\textcolor{darkyellow}{In order to explore when the impact of blade aeroelasticity becomes crucial the deflection contour of the rotor blades at the highest and lowest positions on the rotor has been presented in Fig.\ref{blade_deflection}. Comparing the displacement contours at the highest and the lowest points, shown in Fig.\ref{defl 1} and Fig.\ref{defl 2} respectively, it is observed that the blade undergoes the highest amount of displacement at its lowest point, while at the top the deflection appears to be less considerable. Therefore, moving along the rotor perimeter beginning from the highest to the lowest position the blade deflection sees an upward trend implying that the most crucial position where the effect of flexibility may be found at lower altitudes. However, due to the presence of the tower at the lowest altitude, this blade position is prone to causing the most cruciality in the performance as also demonstrated in the previous section where the highest density of the structures was found to be associated with the part of the flow field behind the tower. It should be noted that the variation in blade deflection between the highest and lowest positions stems from the counterbalancing effect of the weight vector and centrifugal force at higher altitudes, which tend to neutralize each other. Conversely, at lower altitudes, these vectors align, resulting in increased blade deflection.}

\begin{figure}[t]
\nolinenumbers
    \centering

    \begin{subfigure}{0.45\textwidth}
    \nolinenumbers
        \includegraphics[width=\linewidth]{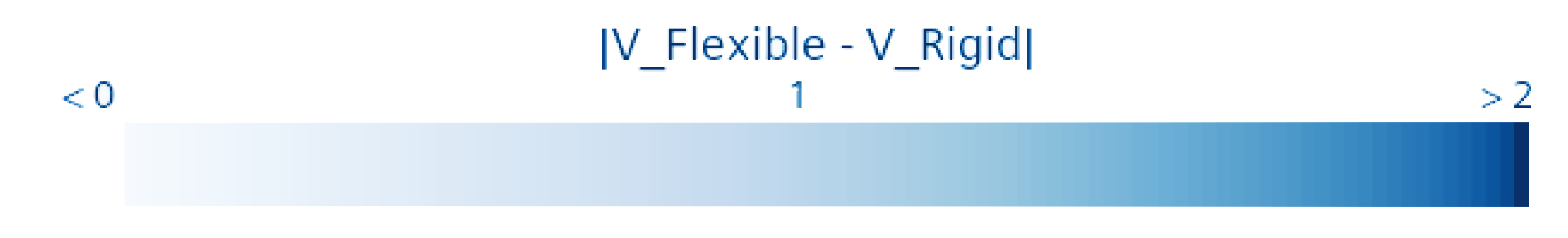} 
    \end{subfigure}
    
    \begin{subfigure}{0.45\textwidth}
    \nolinenumbers
        \includegraphics[width=\linewidth]{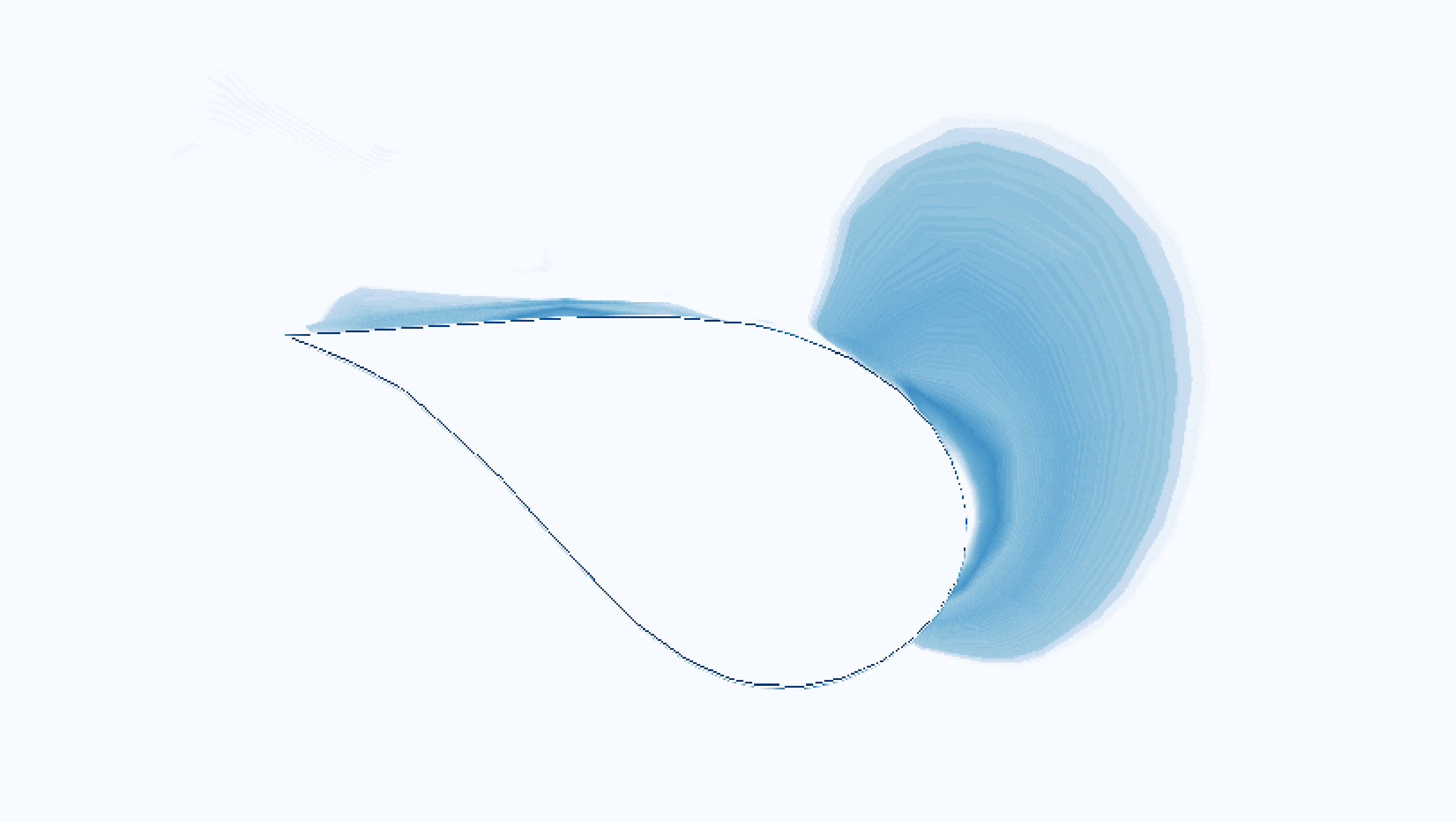}
        \caption{Root Section}
        \label{Root}
    \end{subfigure}
    
    \begin{subfigure}{0.45\textwidth}
    \nolinenumbers
        \includegraphics[width=\linewidth]{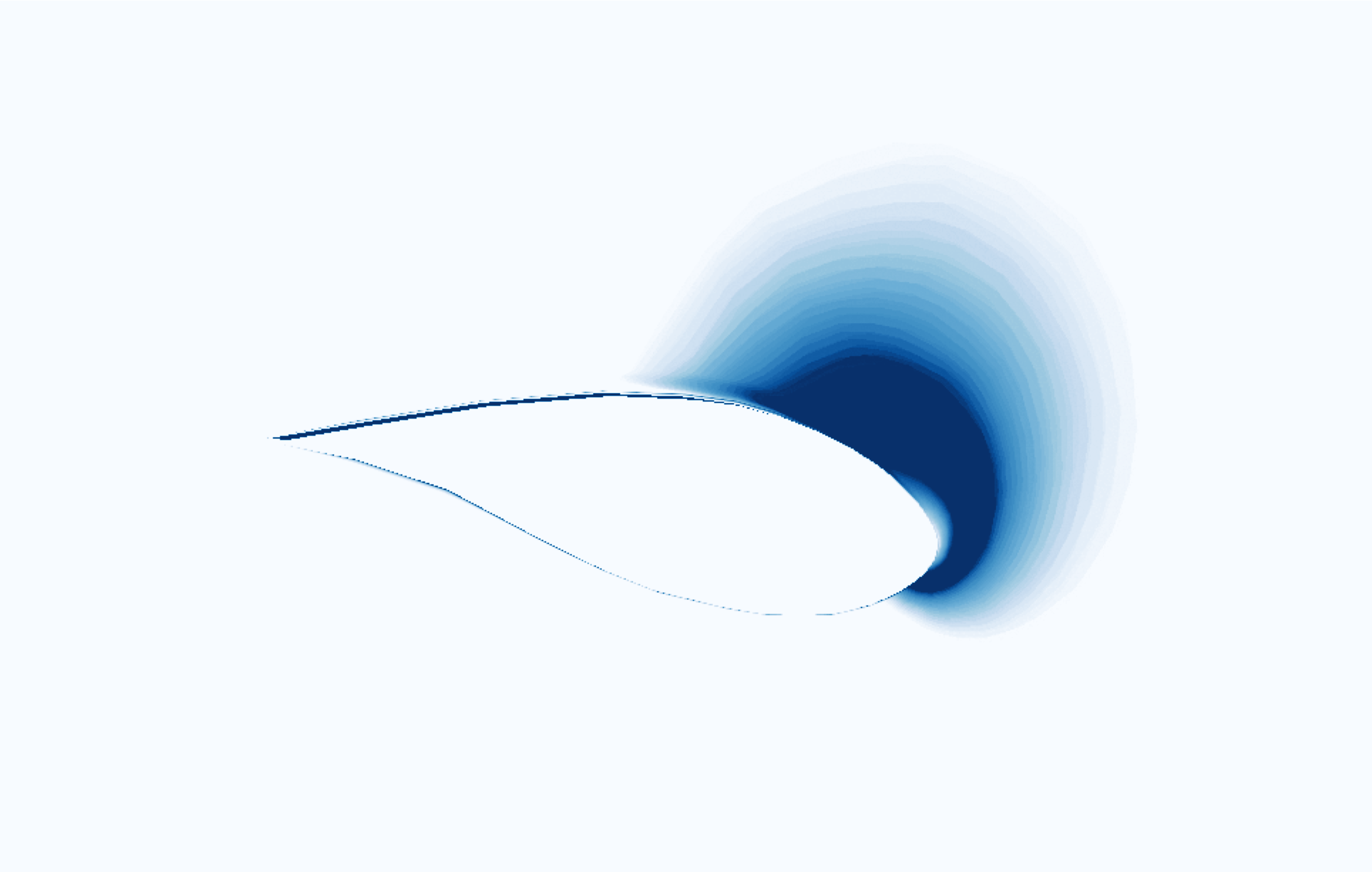} 
        \caption{Middle Section}
        \label{Mid}
    \end{subfigure}

    \begin{subfigure}{0.45\textwidth}
    \nolinenumbers
        \includegraphics[width=\linewidth]{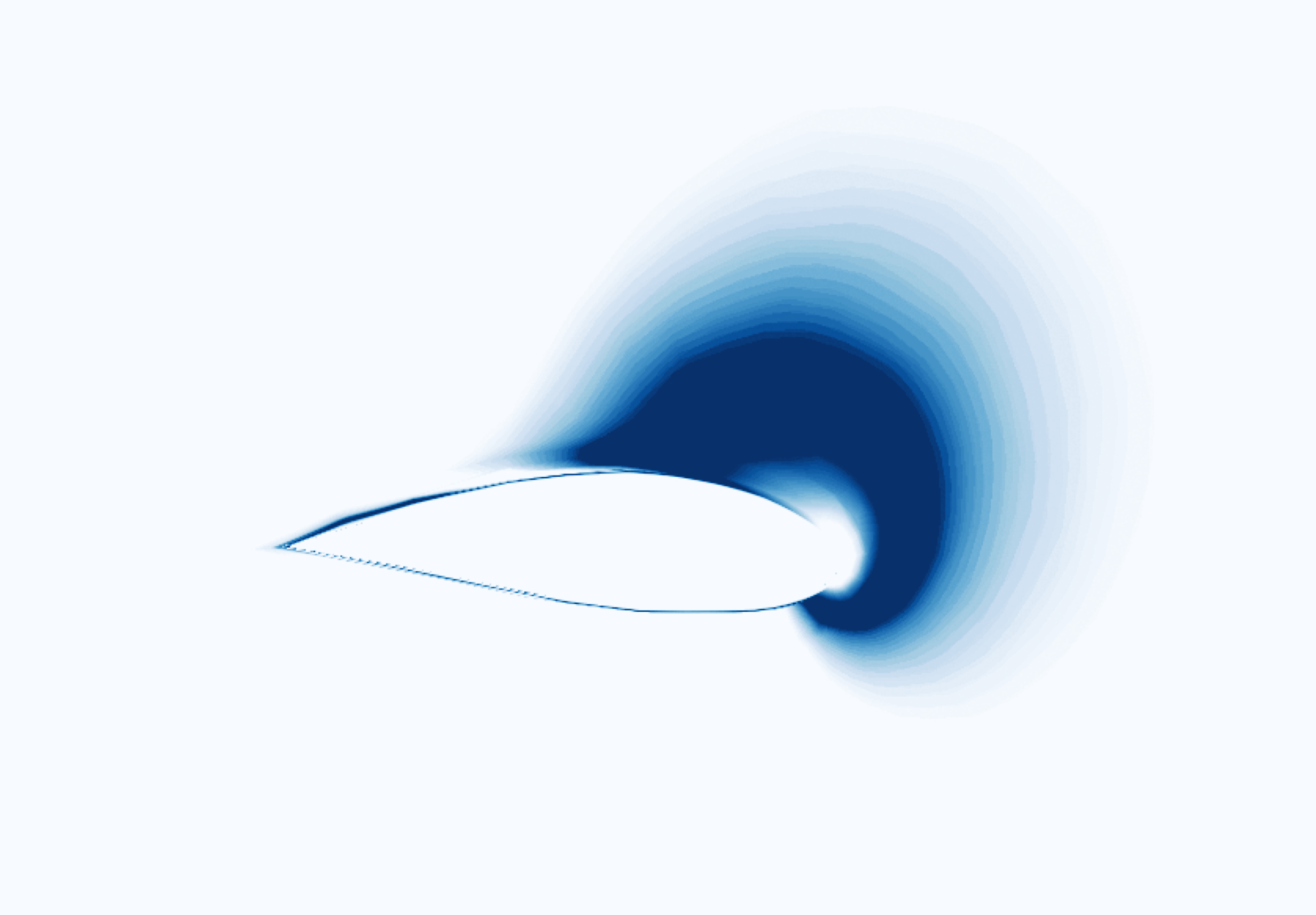} 
        \caption{Tip Section}
        \label{Tip}
    \end{subfigure}

    \caption{\textcolor{darkyellow}{Velocity difference at the root, middle and tip sections}}
    \label{velocity difference}
\end{figure}

\textcolor{darkyellow}{In Fig.\ref{velocity difference} the difference between the velocity fields of the rigid and the flexible cases at the root, mid and tip sections of the blade in the relative frame of motion has been presented. It should be noted that as the trend of the difference at the aforementioned cross sections was almost the same for all the positions of the blade through its course of rotation, the contours were extracted for the blade position where it meets its lowest point on the rotor perimeter. Overall, it is evident in Fig.\ref{velocity difference} that the difference between the flexible and the rigid cases is not highly significant as the maximum difference does not exceed $2(m/s)$ whose effect was also demonstrated in the wake region. First of all, starting from the root section shown in Fig.\ref{Root}, the difference between the flexible and the rigid case seems to be the lowest compared to the other sections of the blade. It is attributed to the fact that the minimum displacement occurs at the root and therefore as expected the velocity field is affected the least in this section. Furthermore, it is observed in Fig.\ref{Root} that the difference is mainly pronounced at the leading edge of the blade which is due to the fact that the blade displacement adds another velocity component to the relative velocity vector causing the relative angle of attack to undergo an alteration proportional to the amount and the direction of the displacement. As a result, the main effect of this mechanism is observed at the leading edge of the blade where the angle of attack plays the main role. In this regard, moving onto the middle and tip section it is observed that analogous to the root section, the leading of both sections demonstrate the highest difference compared to the other parts as seen in Fig.\ref{Mid} and Fig.\ref{Tip} respectively. In addition, the tip section shows the largest velocity difference compared to the other sections due to the highest displacement and the consequent added velocity component. It should be noted that despite the fact that the amount of velocity difference at the middle section is higher than that of the root, it appears from Fig.\ref{Mid} and Fig.\ref{Root} that a more extensive area of the flow field has been affected as result of the blade displacement at the root which stems from the geometry of the airfoil at the root as at the root the airfoil is designed much thicker compared to the other sections which affects a more significant area of the flow field. Therefore, it could be concluded that the blade flexibility mostly affects the near wake region rather than the far wake due to the limited difference in the velocity difference.}

\section{Conclusion}
To analyse the wake behaviour behind a full-scale onshore HAWT a fully coupled two-way FSI simulation has been implemented. The proper orthogonal decomposition has been applied in this matter to extract spatial and temporal modes along with modal energy at the wake region. Results indicate a similar modal behaviour for the rigid and flexible turbines at the wake region which implies the negligible effect of the blade deformation on the wake characteristics. Also, the main source of this deviation is the harmonic nature of the blade deformation through each rotation cycle. 

By observing the wake characteristics at wall-parallel planes behind the turbine, it was concluded that the value of accumulative energy of the first mode (convective mode) decreases as getting further from the turbine. As a consequence, more flow structures in terms of local vortices and fluctuating velocity fields exist at the furthest plane. The flow structures due to the wake shed from the tower tend to move towards the center and merge with that of the nacelle leading to an integral vortical structure 2.5D away from the rotor. In addition, the symmetrical shape of the first mode starts to fade away with the increase in the distance from the rotor proving that while the work extraction is the main characteristic of the planes near the rotor, the effects of wake shed and tip vortices rise into importance in the farther planes. Also, it is concluded that the exclusion of the tower leads to missing a major part of the wake structures, especially at far wake positions.

In summary, the main findings of the paper can be listed as follows:
\begin{itemize}
\item {A similar modal behaviour was observed for the rigid and flexible wind turbines indicating the negligible effect of the blades' flexibility on the wake characteristics behind the HAWT,}
\item {Value of the first mode (convective mode) decreased as getting farther in the wake region,}
\item {The effect of tip vortices caused the elimination of the symmetrical behaviour of the first mode of velocity in locations more than 2D of the rotor,}
\item {The flexibility of the blades does not still affect the wake as much as the tower and nacelle do.}
\item {The modes on the wall-vertical plane demonstrate different behaviours on the near and the far planes.}
\item {It was proven that up to the second the mode the near and far planes show a modal behavior analogous to that of the nacelle and tower on the spanwise plane respectively.}
\item {\textcolor{darkyellow}{It was observed that the highest velocity difference happens at the tip due to the largest blade deflection.}}
\item {\textcolor{darkyellow}{The influence of blade flexibility is primarily concentrated within the near wake region rather than the far wake.}}
\end{itemize}
%


\appendix

\section{\textcolor{darkyellow}{Two-way Coupling}}\label{two-way coupling}
\begin{figure}
\nolinenumbers
	\centering
	\includegraphics[scale=0.6]{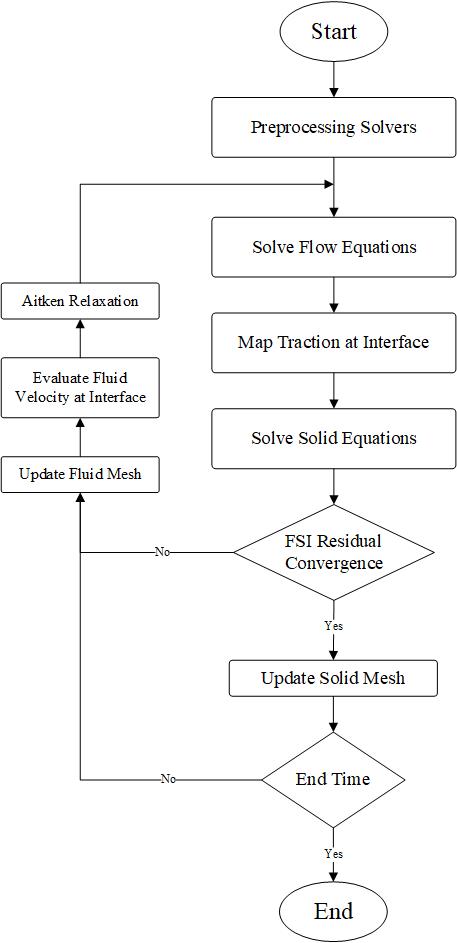}
	\caption{Solution strategy for the partitioned FSI solver}
	\label{FIG_flowchart}
\end{figure}

There exist two methods for numerical simulations of the fluid-structure interaction in terms of dealing with the system of equations. In the so-called \textit{monolithic} approach, the coupled problem is considered as one system of equations, and hence, the flow and structural fields are solved simultaneously. Whereas, in the \textit{partitioned} approach, separate solvers are employed for each of the solid and fluid domains. In the latter, there will be a need for the transfer of forces and displacements between the solvers through time. In fact, there is still no comprehensive discretization method to address both fluid and structural problems accurately. Hence, utilization of the monolithic approach is not common for real-world applications. However, as for the partitioned method, a suitable general-purpose solver can be utilized for each domain in which case we need to implement inter-communication protocols between solvers.\\
In a general two-way FSI problem, after the solid equations are solved, the boundary displacements are imported into the fluid solver and the mesh morpher determines the new locations for the nodes of the fluid domain. Following the flow solution, aerodynamic forces acting on the body are exported to the structural solver. The data of displacement and force can be passed between solvers with \textit{implicit} and \textit{explicit} methods. As for the explicit, data exchange only happens at the end of the time step. While it occurs in the inner loop for the implicit method until convergence takes place. An implicit partitioned approach was utilized for the sake of accuracy and convergence.\\
The main drawback with the conventional partitioned FSI method is that the convergence decreases dramatically as the density ratio of solid to fluid decreases (so-called added mass effects \citep{idelsohn2009fluid, kassiotis2010partitioned, StabileMatthiesBorri2018}, and divergence is observed regardless of the time step value in case of the incompressible flow \citep{fernandez2007projection, van2009added}. To overcome this difficulty, especially in the case of strongly coupled FSI problems, a stabilization method should be employed. \textit{Boundary Interface Added Mass} technique has been incorporated in which the fluid domain traction resulting from the interface deformation and consequently the displaced fluid per unit area are predicted and taken as a force correction in the solver. More details of this technique can be found in the work of Nordenström et al.\citep{nordenstrom2022fluid}.\\
In addition, the starting time is typically regarded as the major bottleneck in the convergence of the simulation of FSI problems with large deflections. To guarantee convergence, almost rigid conditions, i.e. high Young's modules, were applied to the wind turbine at the start of the simulation. Then, the elastic properties decreased to that of the actual values. This corresponds to the oscillatory nature of the results of the generated power depicted in Fig. \ref{Fig_Power}.\\
Finally, to enhance the numerical stability of such an FSI problem with large deformations, Aitken’s dynamical relaxation method is employed \citep{kuttler2008fixed}, in which the values of two previous iterations are used to improve the result of the current iteration to obtain a better convergence. The full procedure of the FSI simulation is shown in Fig. \ref{FIG_flowchart}.

\subsection{\textcolor{darkyellow}{Fluid-Structure Interface}}
The boundary condition at the interface of the solid and the fluid, $\Gamma$, can be defined as follows. Assuming the no-slip condition at the wall, the displacement of the solid, $d_{\Gamma}^S$ is equal to the displacement of the fluid, $d_{\Gamma}^F$, at the interface:
\begin{equation}
\label{eq_solid3}
d_{\Gamma}^S=d_{\Gamma}^F,
\end{equation}
This is true for the rate of change of displacement as well. Simply, substituting the fluid displacement at the boundary with the velocity yields:
\begin{equation}
\label{eq_solid4}
\hat{\mathbf{n}} \cdot \mathbf{v}_{\Gamma}^f=\hat{\mathbf{n}} \cdot \frac{\partial \mathbf{d}_{\Gamma}^s}{\partial t},
\end{equation}
and $\hat{\mathbf{n}}$ is the interface normal vector. Moreover, based on Newton's third law, a balance of forces exists on the interface between the fluid and solid:
\begin{equation}
\label{eq_solid5}
\hat{\mathbf{n}} \cdot \boldsymbol{\sigma}_{\Gamma}^S=\hat{\mathbf{n}} \cdot \boldsymbol{\sigma}_{\Gamma}^F,
\end{equation}
where $\boldsymbol{\sigma}_{\Gamma}^S$ and $\boldsymbol{\sigma}_{\Gamma}^F$ are the solid and fluid stresses at the interface, respectively. In essence, Eq. \ref{eq_solid4} and \ref{eq_solid5} describe the kinematics and dynamics of the fluid-structure interface, respectively.

\section{\textcolor{darkyellow}{Mesh Morpher}}

For the simulation of the strongly coupled FSI problem, a two-way information exchange across the interface is needed. The traction forces across the interface, i.e. the pressure and wall shear, are mapped from the fluid domain to the solid domain. This results in deformations in the solid. The solid deformation across the interface, in turn, is mapped back to the fluid domain. The solid deformation only specifies how the interface nodes are displaced, and hence, the displacements of the inner nodes of the fluid domain must then be determined. To accomplish this, a mesh morphing technique known as \textit{Radial Basis Function} (RBF) has been employed to propagate nodes' displacement throughout the whole fluid domain. Accordingly, a couple of control points are placed on the moving boundary, i.e. fluid-structure interface, at which the displacement is known, and then, the interpolation field will be constructed to specify the displacements of the inner nodes. A full description of the RBF technique can be found in \citep{bao2021vortex}.

\newpage

\bibliography{aipsamp}

\end{document}